\documentclass[10pt]{article}
\usepackage[realmainfile]{currfile}
\usepackage[dvipsnames]{xcolor}
\usepackage{amsmath,amssymb,latexsym}
\usepackage{graphicx}
\usepackage{amscd}
\usepackage{enumitem}
\usepackage{bbm}
\usepackage{bm}
\usepackage[french,english]{babel}
\usepackage{caption}
\usepackage[utf8]{inputenc}
\usepackage[T1]{fontenc}
\usepackage{rotating,amssymb,amsmath,amsfonts,delarray}
\usepackage{wasysym}
\usepackage{geometry}
\usepackage{color}
\usepackage{placeins}
\usepackage{float}
\usepackage{amsthm}

\usepackage{amsfonts}
\usepackage{color}
\usepackage{hyperref}
\usepackage{booktabs}

\usepackage{lipsum} 

\usepackage{graphicx}
\usepackage{caption}
\usepackage{subcaption}

\usepackage{float}

\usepackage{tikz-qtree}
\hypersetup{
	colorlinks=true,
	citecolor=blue,
	linkcolor=blue,
	filecolor=magenta,      
	urlcolor=cyan,
}

\makeatletter
\renewcommand\@biblabel[1]{\textbullet}
\makeatother

\definecolor{amethyst}{rgb}{0.6, 0.4, 0.8}

\newtheorem{rem}{Remark}
\newtheorem{deff}{Definition}
\newtheorem{lem}{Lemma}
\newtheorem{prop}{Proposition}
\usepackage{thmtools}

\newtheorem{theorem}{Theorem}
\newtheorem{cor}{Corollary}

\usepackage[scr=boondoxo,scrscaled=1.05]{mathalfa}

\usepackage{tikz}
\usetikzlibrary{trees}
\usetikzlibrary{shapes}
\usepackage{float}
\usepackage{wrapfig}
\usepackage{xcolor}
\definecolor{amethyst}{rgb}{0.6, 0.4, 0.8}

\definecolor{codegreen}{rgb}{0,0.6,0}
\definecolor{codegray}{rgb}{0.5,0.5,0.5}
\definecolor{codepurple}{rgb}{0.58,0,0.82}
\definecolor{backcolour}{rgb}{0.9,0.9,0.9}

\usepackage{pgfplotstable}

\usepackage[ruled,vlined]{algorithm2e}

\usepackage{multicol}
\usepackage{tikz, pgfplots}
\usepackage{authblk}
\pgfplotsset{compat=1.17}
\usetikzlibrary[arrows.meta,bending, calc, positioning, patterns]
\usetikzlibrary{positioning}

\definecolor{Rcolor}{RGB}{150,160,190}
\newcommand{\Rx}{\fontsize{10pt}{12pt}\selectfont
	\raisebox{.3em}{\hspace{1.2em}%
		\llap{\resizebox{1.09em}{.5em}{\color{black}$\bigcirc$}}%
		\llap{\resizebox{1.199em}{.55em}{\color{darkgray}$\bigcirc$}}%
		\llap{\resizebox{1.19em}{.52em}{\color{gray!50}$\bigcirc$}}%
		\llap{\resizebox{1.1em}{.5em}{\color{gray}$\bigcirc$}}%
		\llap{\resizebox{1.25em}{.55em}{\color{gray}$\bigcirc$}}%
	}%
	\hspace{-.85em}%
	\textbf{%
		\textcolor{black}{\textsf{R}}%
		\hspace{-.025em}\raisebox{.01em}{\llap{\textcolor{Rcolor}{\textsf{R}}}}%
}}%
\newbox\rbox
\savebox\rbox{\scalebox{0.1}{\Rx}}

\makeatletter
\def\R{\scalebox{\f@size}{\usebox\rbox}\xspace}
\makeatletter

\usepackage[onehalfspacing]{setspace}

\begin{document}
	\sloppy
	\title{Tree-structured Ising models under mean parameterization}
	\author{Benjamin Côté$^\dagger$, Hélène Cossette$^*$, Etienne Marceau$^*$ \and 
	\textit{$^\dagger$Department of Statistics and Actuarial Science, University of Waterloo, Ontario, Canada}\\
        \textit{$^*$École d'actuariat, Université Laval, Québec, Canada}
    }
	\maketitle

	\begin{abstract}


    In the risk modeling literature, the Ising model has emerged as a valuable framework for dependent Bernoulli random variables, as its underlying graphical structure captures complex dependence patterns. It is almost always employed in its exponential-family parameterization, which comes with some drawbacks such as the marginal distributions not being fixed and an intractable normalizing constant. 
    We study tree-structured Ising models under their mean parameterization, particularly well suited to risk modeling applications, especially as it eliminates these drawbacks. We study tree-structured Ising models under their mean parameterization, which eliminates these drawbacks and is particularly attractive for risk modeling. The mean parameterization also facilitates the characterization of dependence structures and the distribution of aggregate risks. In particular, we derive an analytic expression for the joint probability generating function, yielding an efficient method for computing the aggregate-risk distribution. Similarly, we obtain the ordinary generating function of expected allocations allowing exact computations for risk allocation. 
    The mean parameterization also allows for a stochastic representation of Ising models, yielding a direct sampling algorithm otherwise unavailable under the exponential-family parameterization. 
    We show that Markov random fields with fixed Poisson marginal distributions provide an efficient and accurate approximation for tree-structured Ising models, in the spirit of Poisson approximations.

	\end{abstract}

	\textbf{Keywords:} Multivariate Bernoulli distribution, Dependence tree, Undirected graphical models, Binary Markov random fields, Poisson approximation, Jungle model, Autologistic model, Markov-Bernoulli distribution. 
	
	\bigskip
	
	\newpage

	\section{Introduction}
	\label{sect:Introduction}

The Ising model stems from works in physics, first introduced in \cite{ising1924beitrag}, where it was proposed as a model to describe macroscopic ferromagnetism phenomena from elementary inter-atomic electromagnetic interaction; for an introduction to the Ising model in physics see \cite{okounkov2022ising}, and for historical remarks see \cite{ising2017fate} and \cite{budrikis2024100}.  
 The popularity of the Ising model has expanded well beyond physics, and it has become a most prominent graphical model with multidisciplinary applications. Textbooks on probabilistic graphical models (among which \cite{koller2009probabilistic}, \cite{maathuis2018handbook}, \cite{theodoridis2015machine} and \cite{wainwright2008graphical}) invariably present the model as a classic example of Markov random fields (MRF). 

\bigskip

Credit risk modeling is one of the research areas outside physics in which the Ising model has prompted interest, see in particular the line of research \cite{molins2005long}, \cite{filiz2012graphical}, \cite{molins2016model}, 
\cite{kato2016long} and
\cite{emonti2025negative} where it is employed to capture the dependence between occurrences of defaults.
Its flexibility given by its underlying graphical structure makes it particularly suitable to encode the complex interactions of credit networks; moreover, the Ising model is the entropy-maximizing distribution when marginal distributions and pairwise correlations are specified.
Other applications of the Ising model in risk modeling include, for example,  \cite{denuit2022conditional} for cyber risk and  \cite{navas2020statistical} for occurrences of avalanches. 
These applications of the Ising model fall within the growing use of graphical models in actuarial science and risk modeling, viz.~\cite{oberoi2020graphical}, \cite{boucher2024modeling} and \cite{cossette2026risk}.

\bigskip

The Ising model is almost always employed 
 in its formulation as a member of the exponential family.  This choice is natural in statistical physics, where the canonical parameters have a physical interpretation: in the most elementary one-parameter Ising model, this parameter represents the inverse temperature of the ferromagnetic metal (\cite{kindermann1980markov}). In many other contexts, however, these parameters lack a direct interpretation.  The exponential-family formulation also comes with an intractable normalizing constant, a drawback that extends to Markov random fields more broadly. These considerations motivate the use of alternative parameterizations with clearer statistical interpretation and improved computational tractability. Yet, in the finance and risk modeling applications mentioned above, the Ising model is employed in its exponential form; and thus  these drawbacks must be contended with. An objective of this paper is to showcase an alternative parameterization well suited to such contexts and its usefulness in deriving results relative to risk-modeling questions.

\bigskip


For dependence modeling purposes, it is desirable to use a parameterization that allows the direct specification of marginal distributions and pairwise correlations, while allowing these quantities to vary independently. The parameterization inherited from physics does not provide this feature as varying a single parameter influences all marginal distributions and pairwise correlations, unless one makes the restrictive assumption of "no external field". This assumption bears limited significance for many statistical applications because it imposes symmetric marginal distributions. In its exponential form, the Ising model arises as an MRF conceived through the fixed-conditional paradigm, as fully discussed in \cite{besag1974}; and this is what prevents this property from holding. 
In Chapter~3 of \cite{wainwright2008graphical}, the authors examine an inversion of this paradigm, letting marginal distributions be fixed, and this yields another parameterization for MRFs, referred to as the mean parameterization. 
The paradigm inversion had also been discussed and performed for a lattice MRF in \cite{pickard1980unilateral}; this idea also underlay the design of the MRF in \cite{pickard1977curious}.  The particular case of the Ising model is treated in Example~3.8 of \cite{wainwright2008graphical}. 
The mean parameterization specifies, directly and separately, marginal distributions and pairwise correlations. 
Connections between parameterizations have already been studied in the literature. 
\cite{teugels1990some} established an alternative parameterization of the joint pmf of the multivariate Bernoulli distribution, directly linked to the marginal distributions and joint moments of $\boldsymbol{J}$, and this parameterization is used in \cite{fontana2018representation}. 
Lemma~1 of \cite{lupparelli2024path} provides the mapping from natural to mean parameters; it
 highlights the non-fixedness of marginals and abstruse impacts on correlation when employing the natural parameterization. Theorem~2.1 of \cite{filiz2012graphical} establishes the existence of a unique converse mapping. See Chapter 5 of \cite{kindermann1980markov} for a discussion of the marginal pmfs of tree-structured Ising models under their canonical parameterization.

\bigskip


In this paper, we focus on Ising models defined on trees, for which the mean parameterization admits a simple characterization. Tree structures offer several advantages for dependence modeling: they facilitate sequential conditioning, as highlighted in \cite{cote2025tree}, and admit efficient structure-learning algorithms \cite{chow1968approximating}, see notably \cite{bresler2020learning} and \cite{nikolakakis2021learning} for more recent developments. The intractable normalizing constant is also eliminated.  An application of tree-structured dependence modeling to risk analysis is presented in \cite{cossette2026risk}.
Tree-structured mean-parameterized Ising models have also been extensively studied in information theory, particularly for broadcasting problems (\cite{evans2000broadcasting}, \cite{mossel2001reconstruction}). 
The questions considered in the broadcasting literature differ, however, from those arising in risk modeling, such as the characterization of dependence structures and the distribution of aggregate risks, which are the primary focus of this article. 

\bigskip

To this end, we derive in Section~\ref{sect:pgfIsing} the joint probability generating function (pgf) of the tree-structured Ising model; the use of the mean parameterization makes this derivation possible. In Section~\ref{sect:PGFcomputation}, we exploit this representation to develop efficient computational methods for evaluating the distribution of the sum of the Ising model's components and for computing expected allocations. Expected allocations play a central role in risk-sharing problems; see \cite{blier-wong2022stochastic}, particularly in conditional mean risk-sharing \cite{denuit2012convex} and capital allocation \cite{tasche2007capital}. \cite{denuit2022conditional} study the conditional mean risk-sharing within an individual risk model constructed with an Ising model to describe the occurrence of claims; their approach, however, requires recourse to simulation to be implemented computationally.

\bigskip

 The mean parameterization leads to the emergence of a clique-based stochastic representation, derived in Section~\ref{sect:montecarlo}. Such a stochastic representation is of particular interest for simulation. In the statistical physics literature, due to the intractable normalizing constant of the exponential form parameterization, simulation of Ising models is typically performed using Markov chain Monte Carlo methods; see Chapter~12 of \cite{koller2009probabilistic} and \cite{izenman2021sampling}. By contrast, the stochastic representation derived in Section~\ref{sect:montecarlo} yields a direct sampling algorithm, and avoids the need for MCMC procedures.
In Section~\ref{sect:PoissonApprox}, we examine how the family of MRFs introduced in \cite{cote2025tree} may approximate, based on the Poisson approximation, tree-structured Ising models. This offers an investigation of the limit behavior of Ising models when modeling very rare events. 
The contributions of Sections~\ref{sect:montecarlo} and~\ref{sect:PoissonApprox} extend beyond computational considerations; they provide further insight into the stochastic dynamics underlying tree-structured Ising models.


\bigskip

The paper is constructed as follows. 
In Section~\ref{sect:parameterizations}, we present and connect the different parameterizations of Ising models and discuss the advantages and drawbacks of each. The next sections present our main results discussed above: in Section~\ref{sect:pgfIsing}, we derive the joint pgf of the Ising model and present efficient computation methods for aggregate risk and risk allocation leveraging this result in Section~\ref{sect:PGFcomputation}. In Section~\ref{sect:montecarlo}, we discuss the clique-based stochastic representation and provide a direct sampling method. Results on the Poisson approximation are provided in Section~\ref{sect:PoissonApprox}. Section~\ref{sect:NumerialExamples} presents a numerical example. Section~\ref{sect:conclusion} concludes.

	\section{Parameterizations}
    \label{sect:parameterizations}

Let $d\in\mathbb{N}_1=\mathbb{N}\backslash\{0\}$ and $\mathcal{V} = \{1,\ldots,d\}$. Define $\mathscr{V}$ as the power set of $\mathcal{V}$, excluding the empty set $\varnothing$. 
	Any $d$-dimensional Bernoulli random vector  $\boldsymbol{J}=(J_v,\, v\in\mathcal{V})$ has a multivariate distribution member of the exponential family with joint probability mass function (pmf) expressible as
\begin{equation}
					 p_{\boldsymbol{J}}(\boldsymbol{x}) = \mathrm{exp}\left({\sum_{\mathcal{W}\in\mathscr{V}} \eta_{\mathcal{W}} \prod_{v\in \mathcal{W}}x_v  - A(\boldsymbol{\eta})}\right)  ,\quad \boldsymbol{x}\in\{0,1\}^{d},
					 \label{eq:multivariateBernoulli-natural}
	\end{equation}
	where $\boldsymbol{\eta} = (\eta_{\mathcal{W}}, \mathcal{W}\in\mathscr{V})$ is the vector of natural parameters and $A(\boldsymbol{\eta})$ serves as a normalizing constant; 
	see \cite{dai2013multivariate}.  We call components of $\boldsymbol{\eta}$ associated to elements of $\mathscr{V}$ comprising $i$ elements the natural parameters of the $i$\textit{th order}, $i\in\{1,\ldots,d\}$.

     \bigskip
      
     Ising models are a subclass of multivariate Bernoulli random vectors. In the upcoming subsections,
     we discuss four different parameterizations of Ising models and how they are intertwined. This allows a better understanding of the assumptions underlying the model and highlights the advantages of each one. Considering the Ising model from different perspectives reveals and connects many of its distributional properties.
     Some notation surrounding graphs and trees are required beforehand.
	
\bigskip
A graph $\mathcal{G}$ is given by the pairing $(\mathcal{V},\mathcal{E})$, where $\mathcal{V}$ represents a set of vertices and $\mathcal{E} \subseteq \mathcal{V}\times\mathcal{V}$ a set of unordered edges.  The neighbors of a vertex, $\mathrm{nei}(u)$, $u\in\mathcal{V}$, are the set of vertices connected to it by an edge: $\mathrm{nei}(u)=\{j\in\mathcal{V}: (u,j)\in\mathcal{E}\}$.  A path from a vertex, $u$, to another, $v$, is a set of edges such that $u$ and $v$ each participate in one edge only while every other participating vertices do so twice or more, and we denote it by $\mathrm{path}(u,v)$.   
	We say that participating vertices are \textit{on} the path. 
	A tree $\mathcal{T}=(\mathcal{V},\mathcal{E})$ is a graph with the specific property that one and only one path exists between every pair of vertices. 
	Labeling one vertex, say $r$, as the root gives the rooted tree $\mathcal{T}_r$. When a tree is rooted, the following filiation relations between vertices become well defined, according to that root:
	 the \textit{descendants} of $u$, $u\in\mathcal{V}$, are the set of vertices such that $u$ is on their path to the root, $\mathrm{dsc}(u) = \{j\in\mathcal{V} : (v, w) \in\mathrm{path}(r,j),\, w\in\mathcal{V} \}$; the \textit{children} of $u$ are the set of vertices that are both its descendant and its neighbor, $\mathrm{ch}(u) = \{j\in\mathcal{V}:j\in \mathrm{nei}(u)\cap \mathrm{dsc}(u)\}$; the parent of $u$ is its sole neighboring vertex that is not its children, $\mathrm{pa}(u)\in\mathrm{nei}(u)\backslash\mathrm{ch}(u)$. The root has no parent since all other vertices are its descendants. 
	One may refer to Chapter 3.3 of \cite{saoub2021graph} for further insight on rooted trees.

\subsection{Natural parameterization}

The Ising model is a special case of the multivariate Bernoulli distribution where higher-order interactions between random variables are dismissed. 
	
		\begin{deff}[Ising model]
		\label{def:Ising}
		A multivariate Bernoulli distribution corresponds to that of an Ising model if its natural parameters of third order and more are all zero. An Ising model is defined with respect to a graph $\mathcal{G}=(\mathcal{V},\mathcal{E})$, where non-zero second order natural parameters dictate its edges,
		that is,  $(u,v) \in\mathcal{E}$ when $\eta_{\{u,v\}} \neq 0$.  An Ising model is tree-structured if $\mathcal{G}$ is a tree. 
	\end{deff}

The probability mass function of an Ising model $\boldsymbol{J}=(J_v,\, v\in\mathcal{V})$ defined on graph $\mathcal{G}=(\mathcal{V},\mathcal{E})$ and with vector of parameters $\boldsymbol{\eta} = (\eta_{\mathcal{W}} ,\; \mathcal{W}\in \mathcal{V}\cup \mathcal{E}) \in \mathbb{R}^{2d-1}$ is given by 
	\begin{equation}
	p_{\boldsymbol{J}}(\boldsymbol{x}) = \mathrm{exp}\left({\sum\limits_{v\in\mathcal{V}} \eta_{v} x_v+ \sum\limits_{(u,v)\in\mathcal{E}} \eta_{(u,v)} x_u x_v   - A(\boldsymbol{\eta})}\right) ,\quad\boldsymbol{x}\in\{0,1\}^d. 
	\label{eq:NaturalIsing}
\end{equation}
We will henceforth refer to (\ref{eq:NaturalIsing}) as the joint pmf in its natural formulation, and  $\boldsymbol{\eta}$ are thus the natural parameters as in (\ref{eq:multivariateBernoulli-natural}). 
Due to 
its multidisciplinary reach, the Ising model has been referred to under other names, varying upon some branches of the literature. Research articles stemming from the work in \cite{besag1974} will call it the \textit{autologistic} (or \textit{auto-logistic}) \textit{model}; articles in statistical science may call it the \textit{binary MRF}, for instance \cite{kaiser2014modeling} and \cite{venema1992modeling}; and it has been employed under the name \textit{Jungle model} in credit risk modeling, starting from \cite{molins2016model}. Although some of the models have slight variations in their formulation, they all refer to the same multivariate Bernoulli distribution.

\bigskip

A consequence of zero higher-order natural parameters is conditional independence properties making the Ising model a MRF. 
	\begin{deff}[MRF]
		\label{def:MRF}
		A vector of random variables $\boldsymbol{X} = (X_v, \, v\in\mathcal{V})$ is a MRF encoded by a graph $\mathcal{G}=(\mathcal{V},\mathcal{E})$ if it satisfies the local Markov property on that graph, meaning
		\begin{equation}
			X_u \perp\!\!\!\perp \left(X_v,\, v\in\mathcal{V}\backslash(\{u\}\cup\mathrm{nei}(u))\right) \mid (X_j,\, j\in\mathrm{nei}(u))
				\label{eq:localMarkovproperty}
		\end{equation}
		 for every vertex $u\in\mathcal{V}$, where $\perp\!\!\!\perp$ marks conditional independence.  We say a MRF is tree-structured if it is encoded by a tree $\mathcal{T}$. 
	\end{deff}


	The Hammersley-Clifford theorem states that a vector of random variables, for which every joint probability is non-zero, is a MRF encrypted on a graph $\mathcal{G}$ if and only if its distribution is Gibbs with respect to $\mathcal{G}$. 
	\begin{deff}[Gibbs distribution]
    \label{def:Gibbs}
		A distribution is Gibbs with respect to a graph $\mathcal{G}=(\mathcal{V},\mathcal{E})$ if it is characterized by a probability density/mass function of the form 
        \begin{equation}
            p(\boldsymbol{x}) = \frac{1}{z}\prod_{\mathcal{W}\in\mathcal{C}}\varphi_{\mathcal{W}}((x_v : v\in \mathcal{W})), \quad \boldsymbol{x}\in\mathbb{R}^d, 
            \label{eq:gibbs}
        \end{equation}
	for some functions $\{\varphi_{\mathcal{W}}: \mathcal{W}\in\mathcal{C}\}$, where $\mathcal{C}$ is the set of cliques of the graph and $z$ is a normalizing constant. Cliques are sets of vertices such that every vertex is a neighbor to every other. 
\end{deff}
Gibbs distributions also satisfy the global Markov property, see Chapter~3 of \cite{lauritzen1996graphical}. On a tree, the global Markov property amounts to the following conditional independencies:
\begin{equation}
X_u \perp\!\!\!\perp X_{v} | \; X_j, \quad \text{with vertex $j$ on }\mathrm{path(u,v)}, \quad \text{for every } u,v\in\mathcal{V}, \,(u,v)\not\in\mathcal{E}. 
    \label{eq:globalMarkov}
\end{equation}

    \begin{lem}
		\label{lem:Ising-Unique}
		All tree-structured MRFs with Bernoulli marginal distributions with non-zero joint probabilities for every configuration on $\{0,1\}^d$ are Ising models. 
	\end{lem}
	\begin{proof}
		The result follows by combining the definition of a Gibbs distribution and that of an Ising model. On a tree, cliques contain one or two vertices only; therefore, for the pmf given in (\ref{eq:multivariateBernoulli-natural}) to factorize over the cliques, the associated Bernoulli multivariate distribution must have null natural parameters of third order or more. 
		Otherwise, such a factorization would be impossible. 
		Hence, it corresponds to the pmf of an Ising model, as given in \eqref{eq:NaturalIsing}. 
	\end{proof}

    Note that the statement would not hold for graphs in general, as they may have higher-dimension cliques allowing the presence of higher-order natural parameters meanwhile satisfying Gibbs factorization. 
  Moreover, Lemma~1 should not be wrongly interpreted as: "All tree-structured graphical models with Bernoulli marginals are Ising models". For example, the binary graphical model proposed in the BDgraph package \cite{mohammadi2019bdgraph}, a model which has wrongly been referred to as an Ising model in \cite{denuit2022conditional}, does not satisfy the conditional independence property \eqref{eq:localMarkovproperty} characterizing a MRF. This model is obtained from a marginal transformation of a Gaussian MRF; this transformation dispels the conditional independence property. It is a member of the family of copula Gaussian graphical models introduced in \cite{dobra2011copula}. Another example is the model discussed in \cite{drton2008binary} and \cite{roverato2013log}; it is not an Ising model as the connected set Markov property does not produce a MRF as per Definition~\ref{def:MRF}.

    \subsection{Canonical parameterization}
    \label{sect:canonical}
In the previous subsection, we defined Ising models on support $\{0,1\}^d$, as the marginals' support $\{0,1\}$  usually represents, in a dependence modeling framework, the occurrence or not of an event at a vertex.
  Ising models, however, more generally encompass all two-state-support models in the form \eqref{eq:NaturalIsing}, whether the two states be two colours (leading to the Potts model as a multicolor extension), signs $\{+,-\}$,  or anything else. In particular, in ferromagnetism applications, components of $\boldsymbol{J}$ represent spins of particles, and take values on $\{-1,1\}$. 
 Obviously, distributional properties are not affected by the support, but this leads to another parameterization for the joint pmf of $\boldsymbol{J}$: 
\begin{equation}
	p_{\boldsymbol{J}}(\boldsymbol{y}) = \mathrm{exp}\left({\sum\limits_{v\in\mathcal{V}} \theta_{v} y_v+ \sum\limits_{(u,v)\in\mathcal{E}} \theta_{(u,v)} y_u y_v   - Z(\boldsymbol{\theta})}\right) ,\quad\boldsymbol{y}\in\{-1,1\}^d, 
    \label{eq:CanonicalIsing-moins1}
\end{equation}
 with $\boldsymbol{\theta}=(\theta_{\mathcal{W}},\;\mathcal{W}\in\mathcal{V}\cup\mathcal{E})\subseteq\mathbb{R}^{2d-1}$. 
This parameterization is usually the one preferred in the physics literature.
The Ising model in \eqref{eq:CanonicalIsing-moins1}, is generally referred to with the qualifier \textit{free boundary} in physics-related works to emphasize that no component of $\boldsymbol{J}$ are degenerate (\cite{kindermann1980markov}, Chapter~5). By Lemma~\ref{lem:Ising-Unique}, there is no ambiguity that the pmf in \eqref{eq:CanonicalIsing-moins1} describes the same distribution as the one in \eqref{eq:NaturalIsing}, with some equivalence between parameters $\boldsymbol{\eta}$ and $\boldsymbol{\theta}$. 


\subsection{Centered parameterization}

In \cite{caragea2009autologistic}, the authors suggest transforming the vertex-idiosyncratic parameters from the natural parameterization into a vector $\boldsymbol{\kappa}=(\kappa_v, \, v\in\mathcal{V})\in(0,1)^d$, defined so that
\begin{equation*}
\eta_v = \ln\left(\frac{\kappa_v}{1-\kappa_v}\right) - \sum_{j\in\mathrm{nei}(v)}\kappa_j \eta_{(v,j)} ,\quad  v\in\mathcal{V}.
\end{equation*}
 This is referred to as the centered parameterization of \eqref{eq:NaturalIsing}, which becomes
\begin{equation}
p_{\boldsymbol{J}}(\boldsymbol{x}) = \mathrm{exp}\left(\sum_{v\in\mathcal{V}}  x_v \ln\left(\frac{\kappa_v}{1-\kappa_v}\right) + \sum_{(u,v)\in\mathcal{E}}\eta_{(u,v)} (x_u- \kappa_u) (x_v-\kappa_v) - B(\boldsymbol{\eta}, \boldsymbol{\kappa})\right),    \quad\boldsymbol{x}\in\{0,1\}^d, 
    \label{eq:CenteredIsing}
\end{equation}
with $\boldsymbol{\kappa}\in(0,1)^d$ and $(\eta_e ,\; e\in\mathcal{E})\in\mathbb{R}^{d-1}$, and where $B(\boldsymbol{\eta}, \boldsymbol{\kappa})$ is a normalizing constant. \cite{caragea2009autologistic} also accounted for covariates, in an autologistic regression perspective, but this is not the focus of the present paper. The advantage of the centered parameterization is that the vector $\boldsymbol{\kappa}$ becomes an approximation of the marginal distributions' means when the interactions given by $(\eta_e ,\; e\in\mathcal{E})$ are not too strong. 
The centered parameterization is discussed for general MRFs in \cite{kaiser2012centered}. 
    
	\subsection{Mean parameterization}
    \label{sect:meanparameterization}

The Ising model is commonly referred to and employed in its canonical parameterization \eqref{eq:CanonicalIsing-moins1}, as the parameters then have a physical interpretation in the context of ferromagnetism. Given its close similarity, the natural parameterization is also often used when modeling dependence between binary random variables, such as in \cite{denuit2022conditional} and \cite{lindberg2016markov}. The centered parameterization has also been employed, for example in \cite{kaiser2014modeling} and \cite{wakeland2017exploring}.
	However, these parameterizations come with some major drawbacks for dependence modeling purposes. First, the
	normalization constants $A(\boldsymbol{\eta})$, $Z(\boldsymbol{\theta})$ and $B(\boldsymbol{\eta}, \boldsymbol{\kappa})$ are arduously computable as they involve a summation over every combination on $\{0,1\}^d$. 
	Second, for the natural formulation, the marginal distributions are not fixed: an alteration of the dependence between $u$ and $v$ through $\eta_{(u,v)}$  results in a change of  
the marginal distributions associated to each vertex in the model.  This is also true for the canonical parameterization, albeit for some special cases (see Remark~\ref{rem:no-external} below), and the centered parameterization, wherein the goodness of the approximation given by $\boldsymbol{\kappa}$ is influenced by $(\eta_e,\,e\in\mathcal{E})$. It has lead to the development in many algorithms to compute likelihoods or marginal probabilities; see \cite{wainwright2008graphical} and \cite{theodoridis2015machine} for an exposition.  

\bigskip

To circumvent these downsides,   \cite{wainwright2008graphical} suggests defining the Ising model instead with 
parameters consisting of marginal and pairwise marginal probabilities, meaning a vector $\boldsymbol{q}=(q_v,\, v\in\mathcal{V})$ such that $q_v=\mathbb{E}[J_v]$ for every $v\in\mathcal{V}$ and a vector $\boldsymbol{\alpha}^{\dagger}=(\alpha_e^{\dagger},\, e\in\mathcal{V})$ such that $\alpha^{\dagger}_{(u,v)}=\mathrm{Cov}(J_u,J_v)$. We rather adopt $\boldsymbol{\alpha}=(\alpha_e,\,e\in\mathcal{E})$  where $\boldsymbol{\alpha}_{(u,v)} = \mathrm{Corr}(J_u,J_v)$ as the vector of pairwise interaction parameters for convenience later on;\footnote{$\mathrm{Corr}$ denotes Pearson's correlation coefficient, $\mathrm{Corr}(J_u,J_v) = (\mathrm{E}[J_uJ_v]-\mathrm{E}[J_u]\mathrm{E}[J_v])/\sqrt{\mathrm{Var}(J_u)\mathrm{Var}(J_v)}$.} $\boldsymbol{q}$ and $\boldsymbol{\alpha}$ are called the mean parameters. 

\bigskip

    The mean parameterization permits the direct derivation of the joint pmf for tree-structured Ising models.
	Rooting the tree yields an intrinsic order of conditioning   according to the filiation relation it defines, making it possible to leverage the conditional independence property \eqref{eq:globalMarkov} for simplification in deriving an expression of the joint pmf. This is made possible from the absence of cycles in a tree. The expression of the joint pmf is given by

    \begin{equation}
    \label{eq:MarkovProba}
        p_{\boldsymbol{J}}(\boldsymbol{x}) = p_{J_r}(x_r) \prod_{v\in\mathrm{dsc}(r)} p_{J_v|J_{\mathrm{pa}(v)} = x_{\mathrm{pa}(v)}}(x_v)  =  p_{J_r}(x_r) \prod_{v\in \mathrm{dsc}(r)} \frac{p_{J_v, J_{\mathrm{pa}(v)}}(x_v, x_{\mathrm{pa}(v)})  }{p_{J_{\mathrm{pa}(v)}}(x_{\mathrm{pa}(v)})}. 
    \end{equation}
	
	\bigskip

From the mean parameterization, marginal pmfs $p_{J_v}$, for $v\in\mathcal{V}$, are thus expressed as $p_{J_v}(x_v) = q_v^{x_v}(1-q_v)^{1-x_v}$, and 
\begin{equation}
    p_{J_v,J_{\mathrm{pa}(v)}}(x_v, x_{\mathrm{pa}(v)}) = p_{J_v}(x_v)p_{J_{\mathrm{pa}(v)}}(x_{\mathrm{pa}(v)}) + \alpha_{(\mathrm{pa}(v),v)} (-1)^{x_v+ x_{\mathrm{pa}(v)}}\sigma_{(\mathrm{pa}(v),v)},
    \label{eq:pair}
\end{equation}
where $\sigma_{(u,v)} = \sqrt{q_uq_v(1-q_u)(1-q_v)}$, for $u,v\in\mathcal{V}$, a notational convention kept throughout the paper. 
Therefore, (\ref{eq:MarkovProba}) becomes
\begin{align}
	p_{\boldsymbol{J}}(\boldsymbol{x})
    &=p_{J_r}(x_r) \prod_{v\in \mathrm{dsc}(r)} \frac{  p_{J_v}(x_v)p_{J_{\mathrm{pa}(v)}}(x_{\mathrm{pa}(v)}) + \alpha_{(\mathrm{pa}(v),v)} (-1)^{x_v+ x_{\mathrm{pa}(v)}}\sigma_{(\mathrm{pa}(v),v)}}{p_{J_{\mathrm{pa}(v)}}(x_{\mathrm{pa}(v)})},
	\label{eq:FixedMarginalisationIsing}
\end{align}
 Obviously, one may choose any root $r$ to define the filiation relations in (\ref{eq:FixedMarginalisationIsing}) and will obtain the same expression for the joint pmf, perhaps only different in appearance. The rooting serves only to specify an order of conditioning. 

\bigskip

Define
$\Lambda(\boldsymbol{q})$ the set of admissible $\boldsymbol{\alpha}$ for given $\boldsymbol{q}\in(0,1)^d$.
For every multiplicand in (\ref{eq:FixedMarginalisationIsing}) to be strictly positive, one requires
 \begin{equation}
 \alpha_{(\mathrm{pa}(v),v)}\in\left(-\min\left(\sqrt{\tfrac{q_{\mathrm{pa}(v)}q_v}{(1-q_{\mathrm{pa}(v)})(1-q_v)}},\sqrt{\tfrac{(1-q_{\mathrm{pa}(v)})(1-q_v)}{q_{\mathrm{pa}(v)}q_v}} \right),\min\left(\sqrt{\tfrac{(1-q_{\mathrm{pa}(v)})q_v}{q_{\mathrm{pa}(v)}(1-q_v)}},\sqrt{\tfrac{q_{\mathrm{pa}(v)}(1-q_v)}{(1-q_{\mathrm{pa}(v)})q_v}}\right)\right),
     \label{eq:alphaconstraints}
 \end{equation}
 for every $v\in\mathcal{V}$, thus bounding the possible dependence schemes according to local interactions. From \eqref{eq:alphaconstraints}, we note that pairwise interactions allow for the entire class of bivariate Bernoulli distributions with fixed marginals of parameters $q_v$ and $q_{\mathrm{pa}(v)}$; see Section 7.1.1 of \cite{joe1997}. 
Proposition 4.1 of \cite{wainwright2008graphical} indicates that the constraints on $\boldsymbol{\alpha}$ given by (\ref{eq:alphaconstraints}), ensuring locally consistent distributions, coincide with those of globally realizable distributions; this means that \eqref{eq:alphaconstraints} entirely characterizes $\Lambda(\boldsymbol{q})$. Consequently, the marginal problem on tree $\mathcal{T}$ (see \cite{jha2025random}) is straightforward in the Bernoulli case. These constraints on $\boldsymbol{\alpha}$ form a convex hull and are greatly discussed in Chapter 3 of \cite{wainwright2008graphical}. 
One recognizes in (\ref{eq:FixedMarginalisationIsing}) a joint pmf belonging to a Gibbs distribution factorizing on $\mathcal{T}$, according to Definition~\ref{def:Gibbs}.
	Given the Hammersley-Clifford theorem and Lemma~\ref{lem:Ising-Unique}, the pmf in (\ref{eq:FixedMarginalisationIsing}) is that of an Ising model.



\bigskip

The mean parameterization for tree-structured Ising models has the advantage of explicitly stating the marginals and the pairwise correlations on edges. There is also no embarrassing normalization constant. Not only are edge correlations specified, but correlations between any pair of components also follow straightforwardly, as shown in the next proposition.

	\begin{prop}
    \label{prop:corr}
       For a tree-structured Ising model $\boldsymbol{J}=(J_v,\, v\in\mathcal{V})$ defined on tree $\mathcal{T}=(\mathcal{V},\mathcal{E})$, formulated through its mean parameterization with vectors $\boldsymbol{q}=(q_v,\,v\in\mathcal{V})\in(0,1)^d$ and $\boldsymbol{\alpha}=(\alpha_e,\,e\in\mathcal{E})\in\Lambda(\boldsymbol{q})$,
  \begin{equation*}
      \mathrm{Corr}(J_u,J_v)=\prod_{e\in\mathrm{path}(u,v)}\alpha_e,\quad u,v\in\mathcal{V}.
  \end{equation*}
	\end{prop}
\begin{proof}
			The proof is provided in Appendix \ref{sect:proofCorr}.
		\end{proof}

The mean parameterization highlights connections between the Ising model and other models proposed in the literature. 
 For instance, we note that the
Markov Bernoulli distribution, introduced in \cite{edwards1960meaning} and studied in \cite{gani1982probability}, \cite{cossette2003ruin} and \cite{cossette2010discrete}, corresponds to an Ising model on a series tree. Their equivalence is clear from Lemma~\ref{lem:Ising-Unique}. In a similar vein, the tree graphical model with conditional independence studied in \cite{padmanabhan2021tree} is an Ising model, again under its mean parameterization. In information theory, tree-structured Ising models correspond to a broadcasting process with asymmetric noisy channels on a tree; see \cite{mossel2001reconstruction}. 

\bigskip

One may define MRFs according to two different paradigms: specify the conditional distributions when every neighbor is known, or fix the marginal distributions. For multivariate Bernoulli distributions, both lead to the Ising model. Natural parameters are linked to the conditional distributions, see equation (4.8) in \cite{besag1974}; mean parameters are linked to the marginals. 
 The two paradigms do not usually lead to the same MRFs in general, and thus the Ising model is an exception. For example, the auto-Poisson MRF in \cite{besag1974}, defined through conditional distributions, differs from the MRF obtained by fixing Poisson marginals in \cite{cote2025tree}. 

    \bigskip

\begin{rem}
\label{rem:no-external}
A common assumption when employing the canonical parameterization is that of ``no external field'', for which
 $\theta_v = 0 $ for all $v\in\mathcal{V}$. In this special case, the mapping between canonical and mean parameters is direct: $q_v = 1/2$ for all $v\in\mathcal{V}$ and $\alpha_e = \tanh \theta_e$ for all $e\in\mathcal{E}$. For this reason, the ensuing model is also sometimes referred to as the symmetric Ising model or the palindromic Ising model (\cite{marchetti2016palindromic}). 
 This particular case's mapping is a well-known fact on the Ising model and often used without proof. One may consult \cite{nikolakakis2021predictive}, Lemmas~12 and~14, for a complete proof. A mapping between canonical and mean parameters becomes highly arduous as soon as the no-external-field assumption is removed, and one must employ marginalization algorithms to proceed, see \cite{wainwright2008graphical} and \cite{theodoridis2015machine}. 
The no-external-field assumption is very limiting for a dependence modeler: for practical applications, marginal distributions are rarely exactly symmetric.  In the statistical literature on the Ising model, however, this assumption is often made -- and the authors hence either assume explicitly or implicitly symmetric marginal distributions. Examples of such research articles include
 \cite{ravikumar2010high}, \cite{daskalakis2019testing}, \cite{bresler2020learning} and \cite{nikolakakis2021predictive}. 
\end{rem}

 

	\section{Joint probability generating function}
	\label{sect:pgfIsing}

	Let $\mathcal{P}_{\boldsymbol{X}}$ denote the joint pgf of a vector of discrete random variables $\boldsymbol{X} = (X_i, i\in \{1,\ldots,d\})$, defined by
	\begin{equation}
		\mathcal{P}_{\boldsymbol{X}}(\boldsymbol{t}) = \mathrm{E}\left[t_1^{X_1}\cdots t_d^{X_d}\right], \quad \quad \boldsymbol{t}=(t_1,\ldots,t_d) \in [-1,1]^d.
		\label{eq:defPGF}
	\end{equation}
      Working with the mean parameterization of the Ising model enables the derivation of the exact expression of its pgf, which we provide in Theorem~\ref{th:PGFIsing} below.
This is because the mean parameterization directly yields marginal and bivariate conditional probabilities, without recourse to marginalization algorithms.  Also, the cumbersome normalization constants $A(\boldsymbol{\eta})$ and $Z(\boldsymbol{\theta})$ of the natural and canonical parameterizations prevent the derivation of an expression of the joint pgf directly from those. 
	
	\begin{theorem}[Joint pgf]
			\label{th:PGFIsing}	
	Consider $\boldsymbol{J}=(J_v,\, v\in\mathcal{V})$, a tree-structured Ising model, mean-parameterized with vectors of parameters $\boldsymbol{q}=(q_v, \, v\in\mathcal{V})\in(0,1)^d$ and $\boldsymbol{\alpha}=(\alpha_e, \, e\in\mathcal{E})\in\Lambda(\boldsymbol{q})$, defined on tree $\mathcal{T}=(\mathcal{V},\mathcal{E})$. Choose a vertex $r\in\mathcal{V}$ as the root of the tree and let $\sigma_{(u,v)} = \sqrt{q_{u}q_{v}(1-q_{u})(1-q_{v})}$ for any $u,v\in\mathcal{V}$. The joint pgf of $\boldsymbol{J}$ is given by
		
		\begin{equation*}
			\mathcal{P}_{\boldsymbol{J}}(\boldsymbol{t}) = (1-q_r)\prod_{i\in\mathrm{ch}(r)}\zeta_i^{\mathcal{T}_r}(\boldsymbol{t}_{i\mathrm{dsc}(i)}) + q_rt_r\prod_{i\in\mathrm{ch}(r)}\xi_i^{\mathcal{T}_r}(\boldsymbol{t}_{i\mathrm{dsc}(i)}),\quad \boldsymbol{t} \in [-1,1]^d, 
		\end{equation*}
	
		with $\zeta_v^{\mathcal{T}_r}(\boldsymbol{t}_{v\mathrm{dsc}(v)})$ and $\xi_v^{\mathcal{T}_r}(\boldsymbol{t}_{v\mathrm{dsc}(v)})$ recursively defined, for $v \in \mathrm{dsc}(r)$, by
		\begin{align}
			  \zeta_v^{\mathcal{T}_r}\left(\boldsymbol{t}_{v\mathrm{dsc}(v)}\right) &= \left(1-q_v+\frac{\alpha_{(\mathrm{pa}(v),v)} \sigma_{(\mathrm{pa}(v),v)}}{1-q_{\mathrm{pa}(v)}}\right)\prod_{i\in\mathrm{ch}(v)}\zeta_i^{\mathcal{T}_r}\left(\boldsymbol{t}_{i\mathrm{dsc}(i)}\right) \notag\\&\quad\quad\quad \quad+ \left(q_v-\frac{\alpha_{(\mathrm{pa}(v),v)}\sigma_{(\mathrm{pa}(v),v)}}{1-q_{\mathrm{pa}(v)}}\right) t_v \prod_{i\in\mathrm{ch}(v)} \xi_i^{\mathcal{T}_r}\left(\boldsymbol{t}_{i\mathrm{dsc}(i)}\right);\label{eq:ZETA}\\
			  \xi_v^{\mathcal{T}_r}\left(\boldsymbol{t}_{v\mathrm{dsc}(v)}\right) &= \left(1-q_v-\frac{\alpha_{(\mathrm{pa}(v),v)} \sigma_{(\mathrm{pa}(v),v)}}{q_{\mathrm{pa}(v)}}\right)\prod_{i\in\mathrm{ch}(v)}\zeta_i^{\mathcal{T}_r}\left(\boldsymbol{t}_{i\mathrm{dsc}(i)}\right)\notag\\ &\quad\quad\quad\quad+ \left(q_v+\frac{\alpha_{(\mathrm{pa}(v),v)} \sigma_{(\mathrm{pa}(v),v)}}{q_{\mathrm{pa}(v)}}\right)t_v\prod_{i\in\mathrm{ch}(v)} \xi_i^{\mathcal{T}_r}\left(\boldsymbol{t}_{i\mathrm{dsc}(i)}\right),\label{eq:XI}
		\end{align}
	where products over an empty set are equal to $1$. 
	\end{theorem}
	
	\begin{proof}
		Reorganizing the factors in the joint pgf gives 
		\begin{equation}
			\mathcal{P}_{\boldsymbol{J}}(\boldsymbol{t}) =  \mathrm{E}\left[\prod_{v\in\mathcal{V}} t_v^{J_v} \right] = \mathrm{E}\left[t_r^{J_r}\prod_{i\in\mathrm{dsc}(r)} t_i^{J_i} \right] = \mathrm{E}\left[t_r^{J_r} \prod_{i\in\mathrm{ch}(r)} \left(t_i^{J_i} \prod_{j\in\mathrm{dsc}(i)} t_j^{J_j} \right) \right],   \notag 
		\end{equation}
		for $\boldsymbol{t}=[-1,1]^d$. We condition on $J_r$, following a Bernoulli distribution of parameter $q_r$ according to mean parameterization, and obtain
		\begin{equation} 
		\mathcal{P}_{\boldsymbol{J}}(\boldsymbol{t}) = (1-q_r) \mathrm{E}\left[\left.\prod_{i\in\mathrm{ch}(r)} \left(t_i^{J_i} \prod_{j\in\mathrm{dsc}(i)} t_j^{J_j} \right) \right| J_r= 0\right] + q_r t_r \mathrm{E}\left[\left.\prod_{i\in\mathrm{ch}(r)} \left(t_i^{J_i} \prod_{j\in\mathrm{dsc}(i)} t_j^{J_j} \right) \right| J_r= 1\right]. \notag
		 \end{equation}
		From the conditional independencies given by the global Markov property, expectations factorize, yielding
		\begin{equation}
			\mathcal{P}_{\boldsymbol{J}}(\boldsymbol{t}) = (1-q_r) \prod_{i\in\mathrm{ch}(r)}\mathrm{E}\left[\left. t_i^{J_i} \prod_{j\in\mathrm{dsc}(i)} t_j^{J_j}  \right| J_r= 0\right] + q_r t_r \prod_{i\in\mathrm{ch}(r)}\mathrm{E}\left[\left. t_i^{J_i} \prod_{j\in\mathrm{dsc}(i)} t_j^{J_j} \right| J_r= 1\right], \label{eq:proofJoingPGF-1}
		\end{equation}
		for $\boldsymbol{t}=[-1,1]^d$. Next, we find recursive relations to resolve the conditional expectations in (\ref{eq:proofJoingPGF-1}). Since $(J_{\mathrm{pa}(v)}, J_v)$ follows a bivariate Bernoulli distribution with dependence parameter $\alpha_{(\mathrm{pa}(v),v)}$, for every $v\in\mathrm{dsc}(r)$,
		the following conditional probabilities may be expressed as such: 
		\begin{align*}
			&\mathrm{Pr}(J_v = 0 | J_{\mathrm{pa}(v)} = 0)  = \left(1-q_v+\frac{\alpha_{(\mathrm{pa}(v),v)} \sigma_{(\mathrm{pa}(v),v)}}{1-q_{\mathrm{pa}(v)}}\right) ; \;  \mathrm{Pr}(J_v = 1 | J_{\mathrm{pa}(v)} = 0) = \left(q_v-\frac{\alpha_{(\mathrm{pa}(v),v)}\sigma_{(\mathrm{pa}(v),v)}}{1-q_{\mathrm{pa}(v)}}\right); \\
			&\mathrm{Pr}(J_v = 0 | J_{\mathrm{pa}(v)} = 1) =  \left(1-q_v-\frac{\alpha_{(\mathrm{pa}(v),v)} \sigma_{(\mathrm{pa}(v),v)}}{q_{\mathrm{pa}(v)}}\right); \;  \mathrm{Pr}(J_v = 1 | J_{\mathrm{pa}(v)} = 1) =  \left(q_v+\frac{\alpha_{(\mathrm{pa}(v),v)} \sigma_{(\mathrm{pa}(v),v)}}{q_{\mathrm{pa}(v)}}\right).
		\end{align*}
		The following relations are then provided by conditioning on $J_{v}$ and applying the global Markov property once more: 
		\begin{align}
			&\mathrm{E}\left[\left.t_v^{J_v} \prod_{j\in\mathrm{dsc}(v)}t_j^{J_j}\right| J_{\mathrm{pa}(v)} = 0 \right]\notag \\&\quad= \left(1-q_v+\frac{\alpha_{(\mathrm{pa}(v),v)} \sigma_{(\mathrm{pa}(v),v)}}{1-q_{\mathrm{pa}(v)}}\right)\prod_{i\in\mathrm{ch}(v)}	\mathrm{E}\left[\left.t_i^{J_i} \prod_{j\in\mathrm{dsc}(i)}t_j^{J_j}\right| J_{\mathrm{pa}(i)} = 0 \right]  \notag\\&\quad\quad+ \left(q_v-\frac{\alpha_{(\mathrm{pa}(v),v)}\sigma_{(\mathrm{pa}(v),v)}}{1-q_{\mathrm{pa}(v)}}\right) t_v \prod_{i\in\mathrm{ch}(v)} 	\mathrm{E}\left[\left.t_i^{J_i} \prod_{j\in\mathrm{dsc}(i)}t_j^{J_j}\right| J_{\mathrm{pa}(i)} = 1 \right] ;\label{eq:expectationZeta}
          \end{align}
          and
          \begin{align}
			&\mathrm{E}\left[\left.t_v^{J_v} \prod_{j\in\mathrm{dsc}(v)}t_j^{J_j}\right| J_{\mathrm{pa}(v)} = 1 \right]\notag \\&\quad= \left(1-q_v-\frac{\alpha_{(\mathrm{pa}(v),v)} \sigma_{(\mathrm{pa}(v),v)}}{q_{\mathrm{pa}(v)}}\right)\prod_{i\in\mathrm{ch}(v)}	\mathrm{E}\left[\left.t_i^{J_i} \prod_{j\in\mathrm{dsc}(i)}t_j^{J_j}\right| J_{\mathrm{pa}(i)} = 0 \right]  \notag\\&\quad\quad+ \left(q_v + \frac{\alpha_{(\mathrm{pa}(v),v)}\sigma_{(\mathrm{pa}(v),v)}}{q_{\mathrm{pa}(v)}}\right) t_v \prod_{i\in\mathrm{ch}(v)} 	\mathrm{E}\left[\left.t_i^{J_i} \prod_{j\in\mathrm{dsc}(i)}t_j^{J_j}\right| J_{\mathrm{pa}(i)} = 1 \right],\label{eq:expectationXi}
		\end{align}
	 with $v\in\mathrm{dsc}(r)$, $\boldsymbol{t}\in[-1,1]^d$. Comparing (\ref{eq:expectationZeta}) and (\ref{eq:expectationXi}) with (\ref{eq:ZETA}) and (\ref{eq:XI}), noting how both pairs of relations recursively intertwine in the same manner, we have 
	  $\mathrm{E}\left[\left.t_v^{J_v} \prod_{j\in\mathrm{dsc}(v)}t_j^{J_j}\right| J_{\mathrm{pa}(v)} = 0 \right] = \zeta_v^{\mathcal{T}_r}(\boldsymbol{t}_{v\mathrm{dsc}(v)})$ and $\mathrm{E}\left[\left.t_v^{J_v} \prod_{j\in\mathrm{dsc}(v)}t_j^{J_j}\right| J_{\mathrm{pa}(v)} = 1 \right] = \xi_v^{\mathcal{T}_r}(\boldsymbol{t}_{v\mathrm{dsc}(v)})$. Substituting these identities in (\ref{eq:proofJoingPGF-1}) grants the desired result. 
	\end{proof}	

	Just as for (\ref{eq:FixedMarginalisationIsing}), one may choose any root $r$ according to which the joint pgf representation in Theorem~\ref{th:PGFIsing} recursively unfolds and obtain the same end result. To specify a root simply determines the filiation relations dictating the order of conditioning in the proof. 

    \bigskip

The heuristics of Theorem~\ref{th:PGFIsing} are similar to those for classical instances of message-passing algorithms, such as the belief propagation algorithm also referred to as the sum-product algorithm; see Chapters~2.5 and~4 of \cite{wainwright2008graphical}.
At each vertex, distributional information on its descendants is retrieved from its children. This information (a \textit{message}) is then \textit{updated} to include the information specific to that vertex and passed along to its parents. Messages thus flow along the edges of the tree from the leaves up to the root.  In Theorem~\ref{th:PGFIsing}, the messages updates correspond to \eqref{eq:ZETA}--\eqref{eq:XI}. Most notably, and perhaps surprisingly, the sum-and-products form of \eqref{eq:ZETA} and \eqref{eq:XI} is reminiscent of the message updates in the belief propagation algorithm, 
despite serving very different purposes. Belief propagation is typically applied to models in their exponential-family form (for Ising, the natural parameterization in \eqref{eq:NaturalIsing}) to recover marginal distributions at each vertex, whereas Theorem~\ref{th:PGFIsing} aims to obtain the joint pgf of a mean-parameterized model (whose marginals are known). Hence, the messages are contributions to probabilities in belief propagation while, in Theorem~\ref{th:PGFIsing}, they are joint pgfs.  
Moreover, the joint pgf being the same for any chosen root, we only need to gather the information at one vertex, and therefore we require messages flowing in one direction for each edge. 



\bigskip

When marginal distributions are not symmetric, the joint-pgf message passed from a vertex $v\in\mathrm{dsc}(r)$ to its parent is different for the two possible outcomes of the Bernoulli distribution at that vertex. This explains the need for both $\zeta_v^{\mathcal{T}_r}$ and $\xi_v^{\mathcal{T}_r}$, to account for each outcome. 
When viewing the Ising model as a broadcast process (Section~\ref{sect:montecarlo} below), 
this difference is explained by an asymmetry in the mutation probabilities depending on the outcomes. 

 
    \bigskip
    
    One may numerically implement the joint pgf from Theorem \ref{th:PGFIsing} by following Algorithm \ref{algo:PGFIsing} below. For such matters, it is most convenient to employ weighted adjacency matrices to provide essential information about the underlying tree in a digestible format to the algorithm.  For a tree $\mathcal{T}=(\mathcal{V},\mathcal{E})$ and vector of parameters $\boldsymbol{\alpha}=(\alpha_e, \, e\in\mathcal{E})$, the corresponding weighted adjacency matrix is given by
\begin{equation*}
	\boldsymbol{A}= (A_{ij},\; (i, j) \in \mathcal{V}\times\mathcal{V}),\quad A_{ij} = \left\{\begin{array}{ll}
		\alpha_{(i,j)}, &(i,j)\in\mathcal{E};\\
		1,&i=j;\\
		0, &\text{elsewhere}.
	\end{array}
	\right.
\end{equation*}
    Here and in all subsequent algorithms, we always assume that weighted adjacency matrices are constructed in topological order: according to a given root, a vertex's associated row and column index comes after its parent's. The root has index 1. This allows for the algorithm to easily retrieve the filiation relations hence defined, as the first non-zero element of each vertex's associated column has row index being their parent's (except for the root's). 
    Designing Algorithm~\ref{algo:PGFIsing} (and the subsequent algorithms) in terms of an adjacency matrix has the advantage of letting it accommodate any arbitrary tree structure. One can use Algorithm E of \cite{cote2025tree} to change the topological order with respect to which the adjacency matrix is constucted. 
    
    \bigskip

\begin{algorithm}[H]
		\label{algo:PGFIsing}
		\caption{Computing the PGF of $\boldsymbol{J}$.} 
		\KwIn{$\boldsymbol{t} = (t_k)_{k\in\mathcal{V}}$;
			  probability parameters $\boldsymbol{q} = (q_k)_{k\in\{1, \ldots, d \}}$; topologically ordered
			 weighted adjacency matrix $\boldsymbol{A} = (A_{ij})_{i\times j\in \mathcal{V}\times \mathcal{V}}$.}
		\KwOut{PGF of $\boldsymbol{J}$ evaluated at $\boldsymbol{t}$.}
		 Construct the $d\times d$ matrices $\boldsymbol{p}^{00}$, $\boldsymbol{p}^{01}$, $\boldsymbol{p}^{10}$ and $\boldsymbol{p}^{11}$ as follows:
		\begin{align*}
			p^{00}_{ij} = \left((1-q_i)(1-q_j) + A_{ij}\sigma_{(i,j)}\right) \mathbbm{1}_{\{A_{ij}\neq0\}}; \quad 
			&p^{01}_{ij} = \left((1-q_i)q_j - A_{ij}\sigma_{(i,j)}\right) \mathbbm{1}_{\{A_{ij}\neq 0\}};\\
			p^{10}_{ij} = \left(q_i(1-q_j) - A_{ij}\sigma_{(i,j)}\right) \mathbbm{1}_{\{A_{ij}\neq 0\}};  \quad 
			&p^{11}_{ij} = \left(q_iq_j + A_{ij}\sigma_{(i,j)}\right) \mathbbm{1}_{\{A_{ij}\neq 0\}}.
		\end{align*}\\
		 Let $\boldsymbol{\zeta}$ and $\boldsymbol{\xi}$ be $d\times d$ matrices of zeros. \\
		 \For{$j = d,\ldots, 2$}{
			 \For{$i = 1, \ldots, j-1$}{
				 Replace $\zeta_{ij}$ and $\xi_{ij}$ by
				\begin{align*}
					\zeta_{ij} &= \left((1-q_j) \prod_{k>j} (\zeta_{kj}+ \mathbbm{1}_{\{A_{kj}=0\}})\right)\mathbbm{1}_{\{A_{ij}\neq0\}},\quad \xi_{ij} = \left(t_jq_j \prod_{k>j} (\xi_{kj}+\mathbbm{1}_{\{A_{kj}=0\}} ) \right) \mathbbm{1}_{\{A_{ij}\neq0\}}. 
				\end{align*}\\
				 Replace $\zeta_{ji}$ and $\xi_{ji}$ by
				\begin{align*}
					\zeta_{ji} &= \frac{p^{00}_ij}{(1-q_i)(1-q_j)}\zeta_{ij} + \frac{p^{01}_{ij}}{(1-q_i)q_j}\xi_{ij},\quad \xi_{ji} = \frac{p^{10}_{ij}}{q_i(1-q_j)}\zeta_{ij} + \frac{p^{11}_{ij}}{q_iq_j}\xi_{ij}.
				\end{align*}\\
			}
		}
		 Compute $\zeta^{(1)} = 	(1-q_1) \prod_{k>1} (\zeta_{k1}+ \mathbbm{1}_{\{A_{k1}=0\}})$.\\
		 Compute $\xi^{(1)} = t_1q_1 \prod_{k>1} (\xi_{k1}+\mathbbm{1}_{\{A_{k1}=0\}} )$.\\
		 Compute $\mathcal{P}_{\boldsymbol{J}} = \zeta^{(1)} + \xi^{(1)}$. \\
		 Return $\mathcal{P}_{\boldsymbol{J}}$.\\
	\end{algorithm}

\section{Pgf-related computation methods}
\label{sect:PGFcomputation}

In this section, we showcase the usefulness, for computational purposes, of having a closed-form expression of the joint pgf. We can efficiently compute the distribution of the sum $K = \sum_{v\in\mathcal{V}} J_v$ and quantities $\mathrm{E}[J_v\mathbbm{1}_{\{K=k\}}]$ for $v\in\mathcal{V}$, $k\in\mathbb{N}$.

	\subsection{Distribution of the sum}

 A naive approach to computing the pmf of the sum, when the joint pmf is known, is to calculate the probabilities for every possible realization of the random vector and regroup the values according to their sum. This approach, however, quickly becomes  computationally expensive, as it requires $2^{d}$  calculations, one for every element of $\{0,1\}^d$. Let us note that this is avoided in the very special case where $\mathcal{T}$ is a star and $\boldsymbol{\alpha} = \alpha\,\boldsymbol{1}_{d-1}$ (where $\boldsymbol{1}_k$ is a vector 1s of length $k$, $k\in\mathbb{N}_1$), for which a closed-form expression of the pmf of $K$ in its natural form is derived in  Section~4.2 of \cite{molins2016model}.

\bigskip

The joint pgf of $\boldsymbol{J}$ allows one to derive the pgf of $K$ through the relation 
\begin{equation*}
    \mathcal{P}_K(t) = \mathcal{P}_{\boldsymbol{J}}(t, \ldots, t), \quad t\in[-1,1].
\end{equation*}
Hence, Theorem~\ref{th:PGFIsing} also provides a recursive expression for the pgf of $K$. Specifically for
     a Markov Bernoulli random vector --- which corresponds to an Ising model defined on a series-tree structure, as stated above in Section~\ref{sect:meanparameterization} ---, the pgf of $K$ is also provided in matrix form in \cite{edwards1960meaning}. The joint pgf given in Theorem~\ref{th:PGFIsing} circumvents the dimensionality issue, while accommodating arbitrary tree structures and parameter values, through Algorithm~\ref{algo:fftM}. It relies on the efficiency of the fast Fourier transform (FFT) algorithm, see \cite{embrechts2009panjer} for more details on this approach.

        \bigskip

\begin{algorithm}[H]
	\label{algo:fftM}
	\caption{Computing the values of the pmf of $K$.} 
	\KwIn{Vector of marginal parameters $\boldsymbol{q}=(q_v)_{v\in\mathcal{V}}$, topologically ordered weighted adjacency matrix $\boldsymbol{A} = (A_{ij})_{i\times j\in \mathcal{V}\times \mathcal{V}}$.}
	\KwOut{Vector $\boldsymbol{p}^{(K)}=(p^{(K)}_k)_{k\in\{0,\ldots,k_{\max}\}}$ such that $p^{(K)}_k$ = $p_{K}(k)$.}
	 Let $n_{\mathrm{FFT}}$ be a power of 2 greater than $d+1$. \\
	 Set $\boldsymbol{b} = (b_{i})_{i\in\{0,\ldots,n_{\mathrm{FFT}}-1\}} = (0,1,0,0,\ldots,0)$.\\
	 Use FFT to compute the discrete Fourier transform $\boldsymbol{\phi}^{(b)}$ of $\boldsymbol{b}$. \\
	 \For{$\ell = 0, 1, \dots, n_{\rm FFT}-1$}{
   Compute $\phi_{\ell}^{(K)} = \mathcal{P}_{J}(\phi_{\ell}^{(b)}, \ldots, \phi_{\ell}^{(b)})$ by employing Algorithm~\ref{algo:PGFIsing}.}
	 Use FFT to compute the inverse discrete Fourier transform $\boldsymbol{p}^{(K)}$ of $\boldsymbol{\phi}^{(K)}$.\\
	 Return $\boldsymbol{p}^{(K)}$.\\
\end{algorithm}

\bigskip

Proposition 3.3 of \cite{padmanabhan2021tree} provides an alternative algorithm for computing the pmf of $K$. Although it does not rely on FFT, it also involves an intertwined recursion. Much like for our approach, the use of the mean parameterization therein plays a key role in enabling the computation of the pmf of $K$.

\subsection{Expected allocations}

 The joint pgf also allows the computation of quantities of the form $\mathbb{E}[J_v\mathbbm{1}_{\{K=k\}}]$, $v\in\mathcal{V}$, $k\in\mathbb{N}$, referred to as expected allocations. As highlighted in \cite{blier2025efficient}, these quantities are key to many computations in risk management, arising for instance in the conditional-mean risk-sharing rule of \cite{denuit2012convex} and in Euler allocations of Tail Value-at-Risk \cite{tasche2007capital}.  See also \cite{denuit2022conditional} in the setting of the Ising model. \cite{blier2025efficient} introduce the ordinary generating function of expected allocations (OGFEA) for efficient computations of $\mathbb{E}[J_v\mathbbm{1}_{\{K=k\}}]$ relying on the FFT algorithm.

\begin{deff}[OGFEA]
Let $\boldsymbol{X} = (X_v,\, v\in\mathcal{V})$ be a vector of random variables taking values in $\mathbb{N}^d$, and let $Y$ denote the sum of its components.
	The ordinary generating function of expected allocations (OGFEA) of $X_v$ with respect to $Y$, $v\in\mathcal{V}$, is $\mathcal{P}_Y^{[v]}$ such that
	\begin{equation*}
		\mathcal{P}_Y^{[v]}(t) = \sum_{k=0}^{\infty}t^k\mathrm{E}[X_v\mathbbm{1}_{\{Y=k\}}], \quad t\in[-1, 1],\,v\in\mathcal{V}.
	\end{equation*}
\end{deff}

We derive the OGFEA of the Ising model from the joint pgf in Theorem~\ref{th:PGFIsing}. Note that here, $\mathrm{E}[J_v\mathbbm{1}_{\{K=k\}}] = \mathrm{Pr}(J_v =1, K=k)$.

\begin{prop}
Let $\boldsymbol{J}=(J_v,\,v\in\mathcal{V})$ be a tree-structured Ising model defined on tree $\mathcal{T}=(\mathcal{V},\mathcal{E})$, mean-parameterized with vectors of parameters $\boldsymbol{q}=(q_v,\,v\in\mathcal{V})\in(0,1)^{d}$ and $\boldsymbol{\alpha}=(\alpha_e,\,e\in\mathcal{E})\in\Lambda(\boldsymbol{q})$.	Assume that the underlying tree is rooted in $j$, $j\in\mathcal{V}$. The OGFEA of $J_j$ with respect to to $K$ is given by 
	\begin{equation*}
		\mathcal{P}_K^{[j]} (t) = q_j t \prod_{i\in\mathrm{ch}(j)} \xi_i^{\mathcal{T}_j}(t_{i\mathrm{dsc}(i)}), \quad t\in[-1,1],
	\end{equation*}
    where $\xi_v^{\mathcal{T}_j}$, for $v\in\mathcal{V}$, is defined as in \eqref{eq:XI}.
\end{prop}
\begin{proof}
	As discussed above, the choice of the root bears no stochastic significance, and one may select any vertex for this purpose when deriving the joint pgf as in Theorem~\ref{th:PGFIsing} without it affecting the result. Choosing $j$ as the root, we notice that $t_j$ is not in $\boldsymbol{t}_{i\mathrm{dsc}(i)}$ for all $i\in\mathrm{ch}(j)$. Therefore, the quantities $\prod_{i\in\mathrm{ch}(j)}\zeta_i^{\mathcal{T}_j}(\boldsymbol{t}_{i\mathrm{dsc}(i)})$ and  $\prod_{i\in\mathrm{ch}(j)}\xi_i^{\mathcal{T}_j}(\boldsymbol{t}_{i\mathrm{dsc}(i)})$ are both simply coefficients of $t_j$. The result then follows from Theorem~2.4 of \cite{blier2025efficient}, stating that  \begin{equation*}
		\mathcal{P}_K^{[v]}(t) = \left. \left(t_j \times \frac{\partial}{\partial t_j} \mathcal{P}_{\boldsymbol{J}}(\boldsymbol{t})\right) \right|_{\boldsymbol{t}=t\,\boldsymbol{1}_d}, \quad t\in[-1,1].  
	\end{equation*}
\end{proof}

The following algorithm complements Algorithm 2 by showing how to leverage the FFT algorithm to compute the expected allocations for the Ising model.\\

\begin{algorithm}[H]
	\label{algo:ogfea}
	\caption{Computing expected allocations $\mathbb{E}[J_v\,\mathbbm{1}_{\{K=k\}}]$, for a chosen $v\in\mathcal{V}$.} 
	\KwIn{Vector of marginal parameters $\boldsymbol{q}=(q_v)_{v\in\mathcal{V}}$, topologically ordered weighted adjacency matrix $\boldsymbol{A} = (A_{ij})_{i\times j\in \mathcal{V}\times \mathcal{V}}$.}
	\KwOut{Vector $\boldsymbol{e}^{(K)}=(e^{(K)}_k)_{k\in\{0,\ldots,d\}}$ such that $e^{(K)}_k$ = $\mathbb{E}[J_v\,\mathbbm{1}_{\{K=k\}}]$.}
    Reorder $\boldsymbol{A}$ so that it is topologically ordered according to vertex $v$. This can be done, for instance, using Algorithm 5 of \cite{cote2025tree}. \\
	  Set $n_{\mathrm{FFT}}$ to be a power of 2 greater than $d+1$. \\
	 Set $\boldsymbol{b} = (b_{i})_{i\in\{0,\ldots,n_{\mathrm{FFT}}-1\}} = (0,1,0,0,\ldots,0)$.\\
	 Use FFT to compute the discrete Fourier transform $\boldsymbol{\phi}^{(b)}$ of $\boldsymbol{b}$. \\
	 \For{$\ell = 0, 1, \dots, n_{\rm FFT}-1$}{
   Compute $\phi_{\ell}^{(K)} = \mathcal{P}_{J}(\phi_{\ell}^{(b)}, \ldots, \phi_{\ell}^{(b)})$ by employing Algorithm~\ref{algo:PGFIsing}, but change the penultimate step by ``Compute $\mathcal{P}_{\boldsymbol{J}}= \xi^{(1)}$.'' }
	 Use FFT to compute the inverse discrete Fourier transform $\boldsymbol{e}^{(K)}$ of $\boldsymbol{\phi}^{(K)}$.\\
	 Return $\boldsymbol{e}^{(K)}$.\\
\end{algorithm}


	\section{Stochastic representation and sampling}
    \label{sect:montecarlo}
	
	In its standard formulations \eqref{eq:NaturalIsing} and \eqref{eq:CanonicalIsing-moins1}, the Ising model does not admit a straightforward sampling procedure because, for each vertex, its dependence relations with every neighbour cannot be decoupled. 
 An implementation of Markov Chain Monte Carlo, called the Gibbs sampler, circumvents this need. It produces a realization of an Ising model
 by iteratively simulating local interactions between random variables, refining the realization's agreement to the joint distribution at each iteration. This, of course, requires many random numbers to produce a single realization of $\boldsymbol{J}$, and may not be efficient when $d$ is large. 
    The Gibbs sampler method is described in detail in Chapter~12 of \cite{koller2009probabilistic}. 
See \cite{izenman2021sampling} for a complete exposition of the sampling problem for an Ising model. 
Most implementations of sampling methods for Ising models employ the Gibbs sampler or variations of it; for instance, the IsingSampler package in R \cite{isingsampler}. 
An alternative to the Gibbs sampler is importance sampling, which consists of weighting realizations of a distribution close to the one from which we aim to sample and for whose a direct sampling method exists. It is more efficient than the Gibbs sampler, but comes with the drawback of accrued variance for Monte Carlo estimations, see Chapter~12 of \cite{koller2009probabilistic}.

    \bigskip

Obviously, a direct sampling method, where one iteratively produces a realization of each component of $\boldsymbol{J}$ using a single random number, would be most efficient. 
This method is made possible for tree-structured Ising models by a 
	 a clique-based stochastic representation, provided in the next theorem. It is enabled by the use of the mean parameterization. Beyond providing a straightforward sampling procedure, the stochastic representation given in Theorem~\ref{th:StoRepr-1}  is a novel result which deepens our understanding of the stochastic dynamics at play in a tree-structured Ising model. 
The stochastic representation illustrates clearly how the change in paradigm discussed in Section~\ref{sect:meanparameterization} affects one's understanding of the model. For a further discussion on how the paradigm change in defining a MRF is reflected in the stochastic representation and sampling procedures, see Section~6 of \cite{cote2025tree}. 

\begin{theorem}
	\label{th:StoRepr-1}
	For a tree $\mathcal{T}=(\mathcal{V},\mathcal{E})$ and a chosen root $r\in\mathcal{V}$, let $\mathcal{T}_r$ be its rooted version. For vectors of parameters $\boldsymbol{q}=(q_v,\, v\in\mathcal{V})\in(0,1)^d$ and $\boldsymbol{\alpha}=(\alpha_e,\, e\in\mathcal{E})\in\Lambda(\boldsymbol{q})$, define $\boldsymbol{J}=(J_v,\, v\in\mathcal{V})$ such that $J_r$ follows a Bernoulli distribution of parameter $q_r$ and, for $v\in\mathrm{dsc}(r)$, we have the construction 
	\begin{equation}
		J_{v} = (1- J_{\mathrm{pa}(v)}) I_{v,0}^{\left\{ q_v - \alpha_{(\mathrm{pa}(v),v)}\frac{\sigma_{(\mathrm{pa}(v),v)}}{1-q_{\mathrm{pa}(v)}}  \right\}} +  J_{\mathrm{pa}(v)}  I_{v,1}^{\left\{ q_v + \alpha_{(\mathrm{pa}(v),v)}\frac{\sigma_{(\mathrm{pa}(v),v)}}{q_{\mathrm{pa}(v)}}  \right\}},
        \label{eq:storep}
	\end{equation}
with $\sigma_{(\mathrm{pa}(v),v)}=\sqrt{q_{\mathrm{pa}(v)} q_v (1-q_{\mathrm{pa}(v)}) (1-q_v)}$ and where $(I_{v,0}^{\{\beta_{v,0}\}},\, v\in\mathcal{V})$ and $(I_{v,1}^{\{\beta_{v,1}\}},\, v\in\mathcal{V})$ are independent vectors whose components are mutually independent and satisfy $I_{v,i}^{\{\beta_{v,i}\}}\sim$~Bernoulli$(\beta_{v,i})$, $v\in\mathcal{V}$, $i\in\{0,1\}$. Then, $\boldsymbol{J}$ is a tree-structured Ising model, designed through its mean parameterization. 
\end{theorem}
\begin{proof}
It is standard to check that $J_v$ is Bernoulli distributed with mean $q_v$ and that $\mathrm{Corr}(J_{\mathrm{pa}(v)},J_v) = \alpha_{(\mathrm{pa}(v),v)}$, for every $v\in\mathcal{V}$. We argue that $\boldsymbol{J}$ is a MRF: For any $v\in\mathcal{V}$, the stochastic representation of $J_v$ depends only on three random variables $J_{\mathrm{pa}(v)}$, $I_{v,0}^{\{\beta_{v,0}\}}$ and $I_{v,1}^{\{\beta_{v,1}\}}$, the latter two being independent of the other idiosyncratic random variables of all other vertex $(I_{w,0}^{\{\beta_{w,0}\}},\, w\in\mathcal{V}\backslash\{v\})$ and $(I_{w,1}^{\{\beta_{w,1}\}},\, w\in\mathcal{V}\backslash\{v\})$. Therefore, dependence relations between random variables are solely induced through the filiation relation $\mathrm{pa}$. For $v\in\mathcal{V}$, dependence is encompassed by its parent random variable and those to whom it is a parent (its children) -- that is, its neighbors -- and this means conditional independence with respect to the neighboring random variables.    
Thus, $\boldsymbol{J}$ satisfies the local Markov property \eqref{eq:localMarkovproperty} on $\mathcal{T}$ and is a MRF according to Definition~\ref{def:MRF}. By Lemma~\ref{lem:Ising-Unique}, it is a tree-structured Ising model. 
\end{proof}

The pmf associated to the stochastic representation in Theorem~\ref{th:StoRepr-1} is, of course, given by \eqref{eq:FixedMarginalisationIsing}. As established before, the choice of the root allows for the filiation relations to be well-defined on the underlying tree, but has no incidence on the resulting distribution for $\boldsymbol{J}$. 

\bigskip

As mentioned above, along with the clique-based stochastic representation in Theorem~\ref{th:StoRepr-1} comes a straightforward simulation algorithm. As for the previous algorithms, the adjacency matrix therein must be constructed in a topological order given a chosen root $r\in\mathcal{V}$. Filiation relations are deduced from such and provide the order of sampling for the components of $\boldsymbol{J}$. Again, that order depends on the chosen root $r$, but that choice has no distributional incidence.

\begin{algorithm}[H]
	\label{algo:simulJi-1}
\caption{Sampling from a tree-structured Ising model.}
\KwIn{probability parameters $\boldsymbol{q}=(q_v)_{v\in\mathcal{V}}$; 
	topologically ordered weighted adjacency matrix $\boldsymbol{A} = (A_{ij})_{i\times j\in \mathcal{V}\times\mathcal{V}}$.}
\KwOut{Realization of $\boldsymbol{J}$.} 
 Simulate $J_1$ as a Bernoulli$(q_1)$ random variable. \\
 \For{$v=2,\ldots,d$}{
	 Compute $\pi_v = \inf\{j:A_{vj}\neq 0\}$.\\
 \If{$J_{\pi_v} = 0$} {
	 Produce a realization of $J_v$ as a $\mathrm{Bernoulli}\left(q_v-A_{\pi_v,v}\tfrac{\sigma_{(\pi_v,v)}}{1-q_{\pi_v}}\right)$ random variable.\\
} 
 \Else{
	 Produce a realization of $J_v$ as a $\mathrm{Bernoulli}\left(q_v+ A_{\pi_v,v} \tfrac{\sigma_{(\pi_v,v)}}{q_{\pi_v}}\right)$ random variable.
}}
	 Return $(J_v,\;v\in\{1,\ldots,d\})$.
\end{algorithm}

\bigskip

Under the fixed-conditional paradigm leading to the other, more common parameterizations, the conception of the stochastic dynamics is not clique-based but neighbor-based. Stochastic dynamics take into account every neighbor's vertex, see \cite{besag1974}; they do not flow down the tree. 
A stochastic representation leaning on that could therefore not produce a straightforward sampling procedure as every vertex's realization depends on the others'; hence the need for the Gibbs sampler. 

\bigskip

As mentioned in Section~\ref{sect:canonical}, the Ising model is more generally conceived as a two-color model, where the color realization of node $v$ is influenced by those of $\mathrm{nei}(v)$. 
 Can we adapt this conception to a stochastic representation flowing down the tree as in \eqref{eq:storep}, but where the stochastic dynamics rather solely consider whether or not a vertex's color is the same as its parent? The answer is no, in general, as asymmetries in the marginal distributions make it necessary to know the precise color of a vertex: it influences the probability of propagation. 
Interestingly, however, it is possible in the special symmetric case where $\boldsymbol{q} = 0.5\,\boldsymbol{1}_d$:
let $J_r \sim $ Bernoulli(0.5) and
\begin{equation}
    J_v = I^{\{\beta_v\}}_v J_{\mathrm{pa}(v)} + \left(1-I^{\{\beta_v\}}_v\right)(1-J_{\mathrm{pa}(v)}), \quad v\in\mathrm{dsc}(r),
    \label{eq:storep-potts}
\end{equation}
where $(I^{\{\beta_v\}}_v,\, v\in\mathrm{dsc}(r))$ is a vector of independent Bernoulli trials of success probability $\beta_v$, and $\beta_v = 0.5(\alpha_{(\mathrm{pa}(v),v)}+1)$ for all $v\in\mathrm{dsc}(r)$. 
That the stochastic representation in \eqref{eq:storep-potts} holds for this specific case may explain its simple mapping between the canonical parameterization and the mean parameterization; this ought to be explored further. If $\alpha_{(\mathrm{pa}(v), v)}>0$ for all $v\in\mathrm{dsc}(r)$, the parameters $\beta_v$ are interpreted as the \textit{mutation probability}, or \textit{flip probability}, in broadcasting problems, see e.g.~\cite{evans2000broadcasting}. 

\bigskip

  As discussed in \cite{fontana2018representation}, all members of the Fréchet class of bivariate symmetric Bernoulli distributions lie on a segment within the convex hull of all bivariate Bernoulli distributions; for an illustration, see Figure~3 of \cite{blier-wong2022stochastic}. Comonotonicity and counter-monotonicity\footnote{Comonotonicity is perfect positive dependence. Counter-monotonicity is perfect negative dependence.} are the extremities of that segment, and independence is its midpoint. This is why \eqref{eq:storep-potts} captures all possible dependence schemes of tree-structured no-external-field Ising models, because it yields a discrete mixture between comonotonicity and counter-monotonicity for all of its pairwise distributions (as in \eqref{eq:gibbs}). In this vein, the stochastic representation described in Section~1.1 of \cite{nam2022ising}, producing a discrete mixture between comonotonicity and independence, allows to describe all of those exhibiting positive dependence. This representation is understood as the broadcasting problem on trees in information theory, which has been of prominent interest, as we mentioned in the Introduction. For edges $e$ with $\alpha_e<0$, one could conversely describe the stochastic dynamics as a mixture between independence and counter-monotonicity.

	\section{Poisson approximation}
	\label{sect:PoissonApprox}
    
Poisson approximation refers to the well-known observation that a Bernoulli distribution and a Poisson distribution
   become increasingly close in total variation distance  as their mean parameter decreases to zero.
It is not only a computational tool for approximating either distribution by the other, but more importantly also a theoretical investigation of the distance between distributions as they tend to limit cases. 
In this section, we extend this observation to the multivariate case for tree-structured Ising models and a family of MRFs with Poisson marginal distributions. 

\bigskip

	Most theoretical works regarding Poisson approximations for dependent Bernoulli trials revolve around the Chen-Stein method from \cite{chen1975poisson}, notably discussed in \cite{arratia1990poisson}; see \cite{barbour1992poisson} and \cite{novak2019poisson} for a review of such works. The Chen-Stein method assesses the error of approximating the sum of dependent Bernoulli random variables by a single Poisson random variable, which could be seen as the sum of independent Poisson random variables. The soundness of the approximation is contingent on the strength of the dependence within the Bernoulli random vector. Multivariate extensions of the Chen-Stein method (e.g.~\cite{pianoforte2023multivariate}, \cite{roos1999rate}) consider the sum of independent Bernoulli random vectors and approximates it by a single Poisson random vector.
    
    \bigskip

    We proceed differently. We look at the similarity between vectors from two families of distribution, one with Bernoulli marginals (the Ising model) and the other with Poisson marginals (the family MPMRF, defined below). 
   In other words, we add dependence between the Poisson random variables in order to mimic the dependence scheme of the Bernoulli trials in the Ising model. We adopt this approach not only to have a multivariate outlook on the relation between Bernoulli and Poisson random variables; our intention is mostly to investigate whether the dependence schemes of these two fixed-marginal MRFs closely relate in terms of their strength and dynamics.

\bigskip

Let $\circ$ denote the binomial thinning operator introduced in \cite{steutel1979discrete}.  For $\alpha\in(0,1)$ and a random variable $X$ taking values in $\mathbb{N}$, binomial thinning is the operation $\alpha \circ X = \sum_{i=1}^X I^{\{\alpha\}}_i$, where $\{I^{\{\alpha\}}_i, \, i\in\mathbb{N}_1\}$ is a sequence of independent Bernoulli trials of success probability $\alpha$. 

\begin{deff}[MPMRF]
\label{def:MPMRF}
    For a tree $\mathcal{T}=(\mathcal{V},\mathcal{E})$, choose a root $r\in\mathcal{V}$ and let $\mathcal{T}_r$ be its rooted version.   Given parameters $\lambda>0$ and $\boldsymbol{\alpha}=(\alpha_e,\,e\in\mathcal{E})\in(0,1)^{d-1}$, let $\boldsymbol{N}=(N_v,\,v\in\mathcal{V})$ be a random vector defined according to the following stochastic representation
\begin{equation}
			N_v = 
            \begin{cases}
                L_r, & \text{if } v=r; \\
                 \alpha_{(\mathrm{pa}(v),v)} \circ N_{\mathrm{pa}(v)} +L_v , & \text{if } v \in dsc(r)
            \end{cases}, \quad v\in\mathcal{V},
            \label{eq:storepMPMRF}
		\end{equation}
        where $\boldsymbol{L}=(L_v,\, v\in\mathcal{V})$ is a vector of independent Poisson random variables, each of mean $\lambda(1-\alpha_{(\mathrm{pa}(v),v)})$, and $\lambda$ for $L_r$. 
The random vector $\boldsymbol{N}$ is a MRF with Poisson marginal distributions of mean $\lambda$.  
    We write $\boldsymbol{N}\sim$ MPMRF$(\lambda,\boldsymbol{\alpha},\mathcal{T})$. 
\end{deff}
By Theorem~2 of \cite{cote2025tree}, the choice of the root has no incidence on the distribution of the vector $\boldsymbol{N}$, hence why we do not specify it for its distributional notation. 
We aim to quantify the proximity between the distributions of $\boldsymbol{N}$ and $\boldsymbol{J}$ when marginal distributions tend to degenerate at 0. This intent is greatly facilitated by the mean parameterization: we directly control the marginal distributions 
via the vector of parameter $\boldsymbol{q}$, thus our observations may be obtained by letting $\boldsymbol{q}$ become small. 

\bigskip

	The total variation distance between two $d$-dimensional vectors of random variables, say $\boldsymbol{X}$ and $\boldsymbol{Y}$, is given by 
	\begin{equation*}
		\mathbbm{d}_{\mathrm{TV}}(\boldsymbol{X},\boldsymbol{Y}) = \sup_{A \in \mathcal{B}(\mathbb{R}^d)}\left\{|\mathrm{Pr}(\boldsymbol{X}\in A) - \mathrm{Pr}(\boldsymbol{Y}\in A) |\right\},
	\end{equation*}
	where $\mathcal{B}(\mathbb{R}^d)$ is the Borel set of $\mathbb{R}^d$. The total variation distance may be employed to measure the closeness of two multivariate distributions. 
	
	\begin{theorem}

		Consider $\boldsymbol{J}$, a tree-structured Ising model defined on tree $\mathcal{T}=(\mathcal{V},\mathcal{E})$, mean-parameterized with vectors of parameters $\boldsymbol{q} = q\,\boldsymbol{1}_{d}$, for $q\in(0,1)$, and $\boldsymbol{\alpha}=(\alpha_e,\, e\in\mathcal{E}) \in (0,1)^{d-1}$. 
		Let $\boldsymbol{N}\sim $ MPMRF$(q,\boldsymbol{\alpha},\mathcal{T})$. As $q$ diminishes, the distribution of $\boldsymbol{N}$ becomes an increasingly accurate approximation of the distribution of $\boldsymbol{J}$, for 
		$\mathbbm{d}_{\rm TV}(\boldsymbol{N},\boldsymbol{J}) \leq dq(1-\mathrm{e}^{-q})$.
		\label{th:PoissonApprox}
	\end{theorem}
	
		\begin{proof}
			The proof is provided in Appendix \ref{sect:proofPoissonApprox}.
		\end{proof}
		
		The approximation given by Theorem \ref{th:PoissonApprox} becomes better as $q$ diminishes not just because the scale of both distributions decreases: the bound obtained on $\mathbbm{d}_{\rm TV}(\boldsymbol{N},\boldsymbol{J})$ is of magnitude $O(q(1-\mathrm{e}^{-q})) = O(q^2)$ and not just $O(q)$. Note that Theorem~\ref{th:PoissonApprox} holds for any vector of dependence parameters $\boldsymbol{\alpha}$ and that the bound for $\mathbbm{d}_{\rm TV}(\boldsymbol{N}, \boldsymbol{J})$ is not affected by its value. The dependence between the elements of $\boldsymbol{J}$ may be strong or not, the distribution of $\boldsymbol{N}$ provides an accurate approximation whenever $q$ is small. 
        
        \bigskip

        As a special case of Theorem~\ref{th:PoissonApprox}, we have that a Markov-Bernoulli distribution may be closely approximated by a Poisson AR(1) distribution as the marginal parameter of each becomes small, given they are the respective distributions of $\boldsymbol{J}$ and $\boldsymbol{N}$ when defined on a series tree. Poisson AR(1) distributions were introduced in \cite{mckenzie1985some}; for a comment on their connection to MPMRF, see Section~2 of \cite{cote2025tree}.  
        A direct consequence of Theorem~\ref{th:PoissonApprox} is that the distributions of the sum of each random vector also grow closer when $q$ becomes small, as made explicit in the following corollary.
		
		
		\begin{cor}
			Under the settings of Theorem \ref{th:PoissonApprox}, as $q$ diminishes, the distribution of $M = \sum_{v\in\mathcal{V}} N_v$ becomes an increasingly accurate approximation of the distribution of $K = \sum_{v\in\mathcal{V}} J_v$, for $\mathbbm{d}_{\mathrm{TV}}(M,K) \leq dq(1-\mathrm{e}^{-q})$.
			\label{th:PoissonApproxSum}
		\end{cor}
		\begin{proof}
			Since the set $\{\boldsymbol{x} \in \mathbb{R}^d :\sum_{v=1}^d x_v \in C\}$, for any $C\in\mathscr{P}(\mathbb{N})$, is an element of $\mathscr{P}(\mathbb{N}^d)$, then 
			\begin{equation*}
				\sup_{C\in\mathscr{P}(\mathbb{N})}\left\{|\mathrm{Pr}(M\in C)-\mathrm{Pr}(K\in C) |\right\} \leq \sup_{A \in \mathscr{P}(\mathbb{N}^d)} \left\{ | \mathrm{Pr}(\boldsymbol{N}\in A) - \mathrm{Pr}(\boldsymbol{J}\in A) | \right\}.
			\end{equation*}
			Hence, $\mathbbm{d}_{\mathrm{TV}}(M,K)\leq dq(1-\mathrm{e}^{-q})$ given Theorem \ref{th:PoissonApprox}. 
		\end{proof}

By Theorem~7 of \cite{cote2025tree}, the random variable $M = \sum_{v\in\mathcal{V}} N_v$ follows a compound Poisson distribution. Compound Poisson distributions provide an alternative to the Poisson distribution for approximating the sum of Bernoulli random variables. Unlike the Poisson distribution, they account for dependence among Bernoulli random variables, making accurate approximations possible even in the presence of strong dependence. This approach is known as a compound Poisson approximation, see \cite{genest2003compound} and the review by \cite{vcekanavivcius2022compound}. Corollary~\ref{th:PoissonApproxSum} establishes such an approximation when the Bernoulli random vector follows an Ising model.   

\bigskip
 
As we show in the next theorem, 
the variability (or riskiness) of $M$ always remains greater than that of $K$, making the Poisson approximation a conservative approximation for a risk modeler. This comparison is made with the convex order, defined as follows: A random variable $Y$ dominates another random variable $X$ with respect to the convex order, noted $X\preceq_{\rm cx}Y$, if $\mathbb{E}[\varphi(X)]\leq\mathbb{E}[\varphi(Y)]$ for all convex functions $\varphi$ such that both expectations are well-defined. In actuarial science and economics, orderings in terms of the convex order mean the embodiment of more risk in distributions, see \cite{rothschild1978increasing}. Textbooks \cite{muller2002} and \cite{shaked2007} offer an introduction to the convex order and to integral stochastic orders in general.

        \begin{theorem}
			Under the settings of Theorem \ref{th:PoissonApprox}, $K\preceq_{cx}M$. 
			\label{th:PoissonApproxConvexOrder}
		\end{theorem}
		\begin{proof}
			The proof is provided in Appendix \ref{sect:proofPoissonApproxConvexOrder}.
		\end{proof}
		

        The Poisson approximation presented in the previous results  deepens our understanding of Ising models by connecting the two families of fixed-marginal MRFs. For instance, as shown by the representation in~\eqref{eq:storepMPMRF}, the distribution of $\boldsymbol{N}$ is predicated on propagation-innovation stochastic dynamics. This comes with nice interpretations for applications as such dynamics are common for network-based phenomena. This behavior does not apply to the Ising model, however, as seen from the stochastic representation in Theorem~\ref{th:StoRepr-1}, but the Poisson approximation shows that a propagation-based interpretation of an Ising model's dynamics would still be valid for rare events on a tree. 
        
        \bigskip

       Another aspect of the theoretical interest of the above results comes in connection with \cite{cote2024centrality}. 
In the latter, the authors analyze topological properties of the family MPMRF. They answer questions such as which component of $\boldsymbol{N}$ has more influence on $M$, depending on the position of its vertex within the tree, or which underlying tree's shape yield a riskier distribution of $M$.    
The proof of the cornerstone result of \cite{cote2024centrality}, Theorem~2, enabling the considerations regarding vertex centrality and tree-shape orderings for the model, relies on the self-decomposability of the Poisson distribution and thus its approach cannot be transposed to the Ising model. The closeness of the joint distributions of the two models given by Theorem~\ref{th:PoissonApprox} signals that these results nonetheless extend to the Ising model, at least when $q$ is small enough.  
        

	\section{Numerical example}
	\label{sect:NumerialExamples}
	
	In this example, we illustrate the results of Sections~\ref{sect:PGFcomputation}, \ref{sect:montecarlo} and~\ref{sect:PoissonApprox}.
    We produce the exact values obtained using the joint pgf and FFT method, relying on Algorithm~\ref{algo:fftM}, and values obtained with Monte Carlo methods, relying on Algorithm~\ref{algo:simulJi-1} to produce samples of $\boldsymbol{J}$. We want to show the relevance of the Poisson approximation: we compute values of the pmf of $K$ and examine how well the pmf of $M$ approximates its values. The values of the pmf of $M$ are obtained by using Algorithm~3 of \cite{cote2025tree}.  
	
	\bigskip
	
	Let $\boldsymbol{J} = (J_1, \ldots, J_7)$ be a vector of random variables following an Ising model with its underlying tree being a binary tree of radius~2. We suppose $\alpha_e = 0.7$ for every $e \in \mathcal{E}$.  
	First, we let $q=0.01$. Results obtained respectively using the exact computation method from Section~\ref{sect:PGFcomputation}, the Poisson approximation and Monte Carlo evaluation, are reported in Table~\ref{tab:ApproxIsing0.01}. As expected given Corollary~\ref{th:PoissonApproxSum}, the values provided by the Poisson approximation are within $dq(1-\mathrm{e}^{-q})=0.00070$ of the exact value, and even suggest that the bound on $\mathbbm{d}_{\mathrm{TV}}(M,K)$ could probably be tighter. 
    We also remark that the Poisson approximation falls within the 90\% confidence intervals for the values obtained from Monte Carlo methods when using $n\in\{1000, 10000\}$ samples. These confidence intervals were derived by repeating the Monte Carlo estimation procedure 1000 times.

	\begin{table}[H]
		\centering
		\begin{tabular}{lcccc}
			\hline
			&$\substack{\text{Exact}\\\text{value}}$&$\substack{\text{Poisson}\\\text{approximation}}$&$\substack{\text{Monte Carlo IC90\%}\\(n=1000)}$&$\substack{\text{Monte Carlo IC90\%}\\(n=10000)}$\\
			\hline	
			$\mathrm{Pr}(K=0)$ &0.97231&0.97239&[0.963, 0.981]&[0.9697, 0.9751]\\
			$\mathrm{Pr}(K=1)$ &0.01309&0.01307&[0.007, 0.019]&[0.0111, 0.0149]\\
			$\mathrm{Pr}(K=2)$ &0.00289&0.00291&[0.000, 0.006]&[0.0021, 0.0037]\\
			$\mathrm{Pr}(K\geq3)$ &0.01170&0.01164&[0.006, 0.017]&[0.0100, 0.0135]\\
			\hline
		\end{tabular}
		\caption{Comparison of the values of the pmf of $K$, the sum of the components of an Ising model with $q = 0.01$, using different computation methods
		} 
		\label{tab:ApproxIsing0.01}
	\end{table}

Consider the function $\pi_X:\mathbb{R}\mapsto \mathbb{R}$ given by $\pi_X(z) = \mathrm{E}[(X-z)_+]$, for $z\in\mathbb{R}$.
For random variables $X$ and $Y$, we have $X\preceq_{\rm cx} Y$ if and only if $\pi_X(z) \leq \pi_Y(z)$ for all $z\in\mathbb{R}$ and $\mathrm{E}[X] = \mathrm{E}[Y]$, see Proposition~3.4.3 of \cite{denuit2006}. 
	One could be interested in the values taken by $\pi_K$ and $\pi_{M}$ to confirm the convex order between $K$ and $M$ as stated in Theorem \ref{th:PoissonApproxConvexOrder}. The values should be very similar, yet those of $M$ should always be superior. Table \ref{tab:ApproxIsing0.01-Stoploss} confirms this. 
	
	\begin{table}[H]
		\centering
		\begin{tabular}{lcc}
			\hline
			&$\substack{\text{Exact}\\\text{value}}$&$\substack{\text{Poisson}\\\text{approximation}}$\\
			\hline	
			$\pi_K(0)$ &0.07000&0.07000\\
			$\pi_K(1)$ &0.04231&0.04239\\
			$\pi_K(2)$ &0.02772&0.02785\\
			$\pi_K(3)$ &0.01602&0.01621\\
			$\pi_K(4)$ &0.00901&0.00922\\
			$\pi_K(5)$ &0.00446&0.00466\\
			$\pi_K(6)$ &0.00121&0.00141\\
			$\pi_K(7)$ &0.00000&0.00016\\
			\hline
		\end{tabular}
		\caption{Values of functions of $\pi_K$ for an Ising model with $q = 0.01$ and its Poisson approximation. 
		} 
		\label{tab:ApproxIsing0.01-Stoploss}
	\end{table}
	
	\bigskip
	
	We choose an even smaller $q$ parameter: $q=0.001$. This should make the Poisson approximation even more performant.  We keep the same underlying tree and the same vector of parameters $\boldsymbol{\alpha} = 0.7\, \boldsymbol{1}_{d-1}$. Computed values of the pmf of $K$ for each of the three methods are presented in Table \ref{tab:ApproxIsing0.001}. Once again, we notice that they are agreeing with Theorem~\ref{th:PoissonApproxSum}, stating that $\mathbbm{d}_{\mathrm{TV}}(M,K)\leq dq(1-\mathrm{e}^{-q})=0.000007$.
	
	\begin{table}[H]
		\centering
		\begin{tabular}{lcccc}
			\toprule
			&$\substack{\text{Exact}\\\text{value}}$&$\substack{\text{Poisson}\\\text{approximation}}$&$\substack{\text{Monte Carlo IC90\%}\\(n=1000)}$&$\substack{\text{Monte Carlo IC90\%}\\(n=10000)}$\\
			\midrule	
			$\mathrm{Pr}(K=0)$ &0.9972031&0.9972039&[0.994, 1.000]&[0.9963, 0.9981]\\
			$\mathrm{Pr}(K=1)$ &0.0013405&0.0013402&[0.000, 0.004]&[0.0008, 0.0019]\\
			$\mathrm{Pr}(K=2)$ &0.0002897&0.0002899&[0.000, 0.001]&[0.0001, 0.0006]\\
			$\mathrm{Pr}(K\geq3)$ &0.0011666&0.0011660&[0.000, 0.003]&[0.0006, 0.0017]\\
			\bottomrule
		\end{tabular}
		\caption{Comparison of the values of the pmf of $K$, the sum of the components of an Ising model with $q = 0.001$, using different computation methods
		} 
		\label{tab:ApproxIsing0.001}
	\end{table}
	
	The similarity between functions $\pi_K$ and $\pi_M$ is also increased, with that of the Poisson approximation being greater only by a jot. This is seen in Table \ref{tab:ApproxIsing0.001-Stoploss}.
	
	\begin{table}[H]
		\centering
		\begin{tabular}{lcc}
			\toprule
			&$\substack{\text{Exact}\\\text{value}}$&$\substack{\text{Poisson}\\\text{approximation}}$\\
			\midrule
			$\pi_K(0)$ &0.007000&0.007000\\
			$\pi_K(1)$ &0.004203&0.004204\\
			$\pi_K(2)$ &0.002747&0.002748\\
			$\pi_K(3)$ &0.001580&0.001582\\
			$\pi_K(4)$ &0.000888&0.000890\\
			$\pi_K(5)$ &0.000438&0.000440\\
			$\pi_K(6)$ &0.000118&0.000120\\
			$\pi_K(7)$ &0.000000&0.000002\\
			\bottomrule
		\end{tabular}
		\caption{Values of functions of $\pi_K$ for an Ising model with $q = 0.001$ and its Poisson approximation.
		} 
		\label{tab:ApproxIsing0.001-Stoploss}
	\end{table}





\section{Conclusion}
\label{sect:conclusion}

The Ising model is a classical century-old graphical model of ubiquitous applications; it has, more recently, been of great interest in risk modeling, notably for credit risk and cyber risk. This paper highlights a lesser known parameterization of the Ising model, the mean parameterization, and exposes the advantages of resorting to this parameterization, notably for such applications. The relevance of this parameterization to risk modeling not only comes from the fact that mean parameters directly indicate marginal distributions and pairwise correlations: expressing the Ising model under its mean parameterization also allows to compute quantities that are particularly relevant to risk modeling. In this vein, we derived the joint pgf of the Ising model under its mean parameterization (in Section~\ref{sect:pgfIsing}), and employed this result to obtain (in Section~\ref{sect:PGFcomputation}) the pgf and pmf of the sum of its components, as well as expected allocations, quantities employed in conditional-mean risk-sharing. The mean parameterization also provides alternative stochastic interpretations of the model, enabling straightforward sampling procedures (Section~\ref{sect:montecarlo}). Finally, the mean parameterization uncovers connections between the Ising model and another family of graphical models employed for risk modeling, the family MPMRF. We proved (in Section~\ref{sect:PoissonApprox}) that one model closely approximates the other when modeling rare events. Such a connection suggests that the centrality and shape-related comparison results from \cite{cote2024centrality} could extend to the Ising model; we leave this as an open question.

	\section*{Acknowledgment}
	This work was partially supported by the Natural Sciences and Engineering Research Council of Canada (Cossette: 04273;  Marceau: 05605; Côté: 581589859) and by the Chaire en actuariat de l'Université Laval (Marceau).
	
	\bibliographystyle{apalike}

	\bibliography{bibPoissonAR1}

\begin{thebibliography}{}

\bibitem[Arratia et~al., 1990]{arratia1990poisson}
Arratia, R., Goldstein, L., and Gordon, L. (1990).
\newblock {P}oisson approximation and the {C}hen-{S}tein method.
\newblock {\em Statistical Science}, pages 403--424.

\bibitem[Barbour et~al., 1992]{barbour1992poisson}
Barbour, A.~D., Holst, L., and Janson, S. (1992).
\newblock {\em Poisson Approximation}.
\newblock Oxford University Press.

\bibitem[Besag, 1974]{besag1974}
Besag, J. (1974).
\newblock Spatial interaction and the statistical analysis of lattice systems.
\newblock {\em Journal of the Royal Statistical Society Series B},
  36(2):192--225.

\bibitem[Blier-Wong et~al., 2022]{blier-wong2022stochastic}
Blier-Wong, C., Cossette, H., and Marceau, E. (2022).
\newblock Stochastic representation of {{FGM}} copulas using multivariate
  {{Bernoulli}} random variables.
\newblock {\em Computational Statistics \& Data Analysis}, 173.

\bibitem[Blier-Wong et~al., 2025]{blier2025efficient}
Blier-Wong, C., Cossette, H., and Marceau, E. (2025).
\newblock Efficient evaluation of risk allocations.
\newblock {\em Insurance: Mathematics and Economics}, 122:119--136.

\bibitem[Boucher et~al., 2024]{boucher2024modeling}
Boucher, J.-P., Crainic, A., LeBlanc, A., and Masse, V. (2024).
\newblock Modeling of fire contagion in farms insurance.
\newblock {\em Variance}, 17(1):1–16.

\bibitem[Bresler and Karzand, 2020]{bresler2020learning}
Bresler, G. and Karzand, M. (2020).
\newblock Learning a tree-structured {Ising} model in order to make
  predictions.
\newblock {\em Annals of Statistics}.

\bibitem[Budrikis, 2024]{budrikis2024100}
Budrikis, Z. (2024).
\newblock 100 years of the {I}sing model.
\newblock {\em Nature Reviews Physics}, 6(9):530--530.

\bibitem[Caragea and Kaiser, 2009]{caragea2009autologistic}
Caragea, P.~C. and Kaiser, M.~S. (2009).
\newblock Autologistic models with interpretable parameters.
\newblock {\em Journal of Agricultural, Biological, and Environmental
  Statistics}, 14:281--300.

\bibitem[{\v{C}}ekanavi{\v{c}}ius and Novak, 2022]{vcekanavivcius2022compound}
{\v{C}}ekanavi{\v{c}}ius, V. and Novak, S. (2022).
\newblock Compound poisson approximation.
\newblock {\em Probability Surveys}, 19:271--350.

\bibitem[Chen, 1975]{chen1975poisson}
Chen, L.~H. (1975).
\newblock Poisson approximation for dependent trials.
\newblock {\em Annals of Probability}, 3(3):534--545.

\bibitem[Chow and Liu, 1968]{chow1968approximating}
Chow, C. and Liu, C. (1968).
\newblock Approximating discrete probability distributions with dependence
  trees.
\newblock {\em IEEE transactions on Information Theory}, 14(3):462--467.

\bibitem[Cossette et~al., 2026]{cossette2026risk}
Cossette, H., C{\^o}t{\'e}, B., Dubeau, A., and Marceau, E. (2026).
\newblock On a risk model with tree-structured {P}oisson {M}arkov random field
  frequency, with application to rainfall events.
\newblock {\em ASTIN Bulletin: The Journal of the IAA}, 56(2):317--343.

\bibitem[Cossette et~al., 2003]{cossette2003ruin}
Cossette, H., Landriault, D., and Marceau, E. (2003).
\newblock Ruin probabilities in the compound {M}arkov binomial model.
\newblock {\em Scandinavian Actuarial Journal}, 2003(4):301--323.

\bibitem[Cossette et~al., 2010]{cossette2010discrete}
Cossette, H., Marceau, E., and Maume-Deschamps, V. (2010).
\newblock Discrete-time risk models based on time series for count random
  variables.
\newblock {\em ASTIN Bulletin: The Journal of the IAA}, 40(1):123--150.

\bibitem[C{\^o}t{\'e} et~al., 2024]{cote2024centrality}
C{\^o}t{\'e}, B., Cossette, H., and Marceau, E. (2024).
\newblock Centrality and shape-related comparisons in a tree-structured
  {M}arkov random field.
\newblock {\em arXiv preprint arXiv:2410.20240}.

\bibitem[C{\^o}t{\'e} et~al., 2025]{cote2025tree}
C{\^o}t{\'e}, B., Cossette, H., and Marceau, E. (2025).
\newblock Tree-structured {M}arkov random fields with {P}oisson marginal
  distributions.
\newblock {\em Journal of Multivariate Analysis}, page 105418.

\bibitem[Dai et~al., 2013]{dai2013multivariate}
Dai, B., Ding, S., Wahba, G., et~al. (2013).
\newblock Multivariate {B}ernoulli distribution.
\newblock {\em Bernoulli}, 19(4):1465--1483.

\bibitem[Daskalakis et~al., 2019]{daskalakis2019testing}
Daskalakis, C., Dikkala, N., and Kamath, G. (2019).
\newblock Testing {I}sing models.
\newblock {\em IEEE Transactions on Information Theory}, 65(11):6829--6852.

\bibitem[Denuit and Dhaene, 2012]{denuit2012convex}
Denuit, M. and Dhaene, J. (2012).
\newblock Convex order and comonotonic conditional mean risk sharing.
\newblock {\em Insurance: Mathematics and Economics}, 51(2):265--270.

\bibitem[Denuit et~al., 2006]{denuit2006}
Denuit, M., Dhaene, J., Goovaerts, M., and Kaas, R. (2006).
\newblock {\em Actuarial Theory for Dependent Risks: Measures, Orders and
  Models}.
\newblock John Wiley \& Sons.

\bibitem[Denuit and Robert, 2022]{denuit2022conditional}
Denuit, M. and Robert, C.~Y. (2022).
\newblock Conditional mean risk sharing in the individual model with graphical
  dependencies.
\newblock {\em Annals of Actuarial Science}, 16(1):183--209.

\bibitem[Dobra and Lenkoski, 2011]{dobra2011copula}
Dobra, A. and Lenkoski, A. (2011).
\newblock Copula gaussian graphical models and their application to modeling
  functional disability data.
\newblock {\em Annals of Applied Statistics}, 5(2):969--993.

\bibitem[Drton and Richardson, 2008]{drton2008binary}
Drton, M. and Richardson, T.~S. (2008).
\newblock Binary models for marginal independence.
\newblock {\em Journal of the Royal Statistical Society Series B},
  70(2):287--309.

\bibitem[Edwards, 1960]{edwards1960meaning}
Edwards, A. (1960).
\newblock The meaning of binomial distribution.
\newblock {\em Nature}, 186(4730):1074--1074.

\bibitem[Embrechts and Frei, 2009]{embrechts2009panjer}
Embrechts, P. and Frei, M. (2009).
\newblock Panjer recursion versus {{FFT}} for compound distributions.
\newblock {\em Mathematical Methods of Operations Research}, 69(3):497--508.

\bibitem[Emonti and Fontana, 2025]{emonti2025negative}
Emonti, C. and Fontana, R. (2025).
\newblock Negative correlations in {I}sing models.
\newblock {\em Statistics for Innovation III: SIS 2025, Short Papers,
  Contributed Sessions 2}, page 182.

\bibitem[Epskamp, 2023]{isingsampler}
Epskamp, S. (2023).
\newblock {\em Ising{S}ampler: Sampling Methods and Distribution Functions for
  the {I}sing Model}.
\newblock R package version 0.2.3.

\bibitem[Evans et~al., 2000]{evans2000broadcasting}
Evans, W., Kenyon, C., Peres, Y., and Schulman, L.~J. (2000).
\newblock Broadcasting on trees and the {I}sing model.
\newblock {\em Annals of Applied Probability}, pages 410--433.

\bibitem[Filiz et~al., 2012]{filiz2012graphical}
Filiz, I.~O., Guo, X., Morton, J., and Sturmfels, B. (2012).
\newblock Graphical models for correlated defaults.
\newblock {\em Mathematical Finance}, 22(4):621--644.

\bibitem[Fontana and Semeraro, 2018]{fontana2018representation}
Fontana, R. and Semeraro, P. (2018).
\newblock Representation of multivariate {B}ernoulli distributions with a given
  set of specified moments.
\newblock {\em Journal of Multivariate Analysis}, 168:290--303.

\bibitem[Gani, 1982]{gani1982probability}
Gani, J. (1982).
\newblock On the {{Probability Generating Function}} of the {{Sum}} of
  {{{M}arkov Bernoulli Random Variables}}.
\newblock {\em Journal of Applied Probability}, 19:321--326.

\bibitem[Genest et~al., 2003]{genest2003compound}
Genest, C., Marceau, E., and Mesfioui, M. (2003).
\newblock Compound {P}oisson approximations for individual models with
  dependent risks.
\newblock {\em Insurance: Mathematics and Economics}, 32(1):73--91.

\bibitem[Ising, 1924]{ising1924beitrag}
Ising, E. (1924).
\newblock {\em Beitrag zur Theorie des Ferro-und Paramagnetismus}.
\newblock PhD thesis, Grefe \& Tiedemann Hamburg, Germany.

\bibitem[Ising et~al., 2017]{ising2017fate}
Ising, T., Folk, R., Kenna, R., Berche, B., and Holovatch, Y. (2017).
\newblock The fate of {E}rnst {I}sing and the fate of his model.
\newblock {\em arXiv preprint arXiv:1706.01764}.

\bibitem[Izenman, 2021]{izenman2021sampling}
Izenman, A.~J. (2021).
\newblock Sampling algorithms for discrete {M}arkov random fields and related
  graphical models.
\newblock {\em Journal of the American Statistical Association},
  116(536):2065--2086.

\bibitem[Jha and Nechita, 2025]{jha2025random}
Jha, A.~K. and Nechita, I. (2025).
\newblock On random classical marginal problems with applications to quantum
  information theory.
\newblock {\em Journal of Physics A: Mathematical and Theoretical},
  58(8):085302.

\bibitem[Joe, 1997]{joe1997}
Joe, H. (1997).
\newblock {\em Multivariate Models and Multivariate Dependence Concepts}.
\newblock CRC Press.

\bibitem[Kaiser et~al., 2012]{kaiser2012centered}
Kaiser, M.~S., Caragea, P.~C., and Furukawa, K. (2012).
\newblock Centered parameterizations and dependence limitations in {M}arkov
  random field models.
\newblock {\em Journal of Statistical Planning and Inference},
  142(7):1855--1863.

\bibitem[Kaiser et~al., 2014]{kaiser2014modeling}
Kaiser, M.~S., Pazdernik, K.~T., Lock, A.~B., and Nutter, F.~W. (2014).
\newblock Modeling the spread of plant disease using a sequence of binary
  random fields with absorbing states.
\newblock {\em Spatial Statistics}, 9:38--50.

\bibitem[Kato, 2016]{kato2016long}
Kato, K. (2016).
\newblock Long-range {I}sing model for credit portfolios with heterogeneous
  credit exposures.
\newblock {\em Physica A: Statistical Mechanics and its Applications},
  462:1103--1119.

\bibitem[Kindermann and Snell, 1980]{kindermann1980markov}
Kindermann, R. and Snell, J.~L. (1980).
\newblock {\em {M}arkov random fields and their applications}, volume~1.
\newblock American Mathematical Society.

\bibitem[Koller and Friedman, 2009]{koller2009probabilistic}
Koller, D. and Friedman, N. (2009).
\newblock {\em {Probabilistic Graphical Models: Principles and Techniques}}.
\newblock MIT press.

\bibitem[Lauritzen, 1996]{lauritzen1996graphical}
Lauritzen, S.~L. (1996).
\newblock {\em Graphical Models}.
\newblock Clarendon Press.

\bibitem[Lindberg, 2016]{lindberg2016markov}
Lindberg, O. (2016).
\newblock Markov random fields in cancer mutation dependencies.
\newblock Master's thesis, University of Turku.

\bibitem[Lupparelli and Marchetti, 2024]{lupparelli2024path}
Lupparelli, M. and Marchetti, G.~M. (2024).
\newblock Path-dependent parametric decompositions in {I}sing models.
\newblock {\em Electronic Journal of Statistics}, 18(2):4382--4406.

\bibitem[Maathuis et~al., 2018]{maathuis2018handbook}
Maathuis, M., Drton, M., Lauritzen, S., and Wainwright, M. (2018).
\newblock {\em Handbook of Graphical Models}.
\newblock CRC Press.

\bibitem[Marchetti and Wermuth, 2016]{marchetti2016palindromic}
Marchetti, G.~M. and Wermuth, N. (2016).
\newblock Palindromic {B}ernoulli distributions.
\newblock {\em Electronic Journal of Statistics}, 10:2435--2460.

\bibitem[McKenzie, 1985]{mckenzie1985some}
McKenzie, E. (1985).
\newblock Some simple models for discrete variate time series 1.
\newblock {\em Journal of the American Water Resources Association},
  21(4):645--650.

\bibitem[Mohammadi and Wit, 2019]{mohammadi2019bdgraph}
Mohammadi, R. and Wit, E.~C. (2019).
\newblock {BD}graph: An {R} package for {B}ayesian structure learning in
  graphical models.
\newblock {\em Journal of Statistical Software}, 89:1--30.

\bibitem[Molins and Vives, 2005]{molins2005long}
Molins, J. and Vives, E. (2005).
\newblock Long range {I}sing model for credit risk modeling.
\newblock In {\em AIP Conference Proceedings}, pages 156--161.

\bibitem[Molins and Vives, 2016]{molins2016model}
Molins, J. and Vives, E. (2016).
\newblock Model risk on credit risk.
\newblock {\em Risk and Decision Analysis}, 6(1):65--78.

\bibitem[Mossel, 2001]{mossel2001reconstruction}
Mossel, E. (2001).
\newblock Reconstruction on trees: beating the second eigenvalue.
\newblock {\em Annals of Applied Probability}, 11(1):285--300.

\bibitem[M{\"u}ller and Stoyan, 2002]{muller2002}
M{\"u}ller, A. and Stoyan, D. (2002).
\newblock {\em Comparison Methods for Stochastic Models and Risks}.
\newblock Wiley.

\bibitem[Nam et~al., 2022]{nam2022ising}
Nam, D., Sly, A., and Zhang, L. (2022).
\newblock Ising model on trees and factors of iid.
\newblock {\em Communications in Mathematical Physics}, 389(2):1009--1046.

\bibitem[Navas~Portella, 2020]{navas2020statistical}
Navas~Portella, V. (2020).
\newblock {\em Statistical Modelling of Avalanche Observables: Criticality and
  Universality}.
\newblock PhD thesis, Universitat de Barcelona.

\bibitem[Nikolakakis, 2021]{nikolakakis2021learning}
Nikolakakis, K. (2021).
\newblock {\em Learning Tree-Structured Models from Noisy Data}.
\newblock PhD thesis, Rutgers The State University of New Jersey, School of
  Graduate Studies.

\bibitem[Nikolakakis et~al., 2021]{nikolakakis2021predictive}
Nikolakakis, K.~E., Kalogerias, D.~S., and Sarwate, A.~D. (2021).
\newblock Predictive learning on hidden tree-structured {I}sing models.
\newblock {\em Journal of Machine Learning Research}, 22(1):2713--2794.

\bibitem[Novak, 2019]{novak2019poisson}
Novak, S. (2019).
\newblock Poisson approximation.
\newblock {\em Probability Surveys}, 16:937--942.

\bibitem[Oberoi et~al., 2020]{oberoi2020graphical}
Oberoi, J.~S., Pittea, A., and Tapadar, P. (2020).
\newblock A graphical model approach to simulating economic variables over long
  horizons.
\newblock {\em Annals of Actuarial Science}, 14(1):20--41.

\bibitem[Okounkov, 2022]{okounkov2022ising}
Okounkov, A. (2022).
\newblock The {I}sing model in our dimension and our times.
\newblock In {\em International Congress of Mathematicians, Prize Lectures},
  volume~1, pages 376--413.

\bibitem[Padmanabhan and Natarajan, 2021]{padmanabhan2021tree}
Padmanabhan, D. and Natarajan, K. (2021).
\newblock Tree bounds for sums of {B}ernoulli random variables: A linear
  optimization approach.
\newblock {\em Informs Journal on Optimization}, 3(1):23--45.

\bibitem[Pianoforte and Turin, 2023]{pianoforte2023multivariate}
Pianoforte, F. and Turin, R. (2023).
\newblock Multivariate {P}oisson and {P}oisson process approximations with
  applications to {B}ernoulli sums and-statistics.
\newblock {\em Journal of Applied Probability}, 60(1):223--240.

\bibitem[Pickard, 1977]{pickard1977curious}
Pickard, D. (1977).
\newblock A curious binary lattice process.
\newblock {\em Journal of Applied Probability}, 14(4):717--731.

\bibitem[Pickard, 1980]{pickard1980unilateral}
Pickard, D.~K. (1980).
\newblock Unilateral {M}arkov fields.
\newblock {\em Advances in Applied Probability}, 12(3):655--671.

\bibitem[Ravikumar et~al., 2010]{ravikumar2010high}
Ravikumar, P., Wainwright, M.~J., and Lafferty, J.~D. (2010).
\newblock High-dimensional {I}sing model selection using $\ell_1$-regularized
  logistic regression.
\newblock {\em Annals of Statistics}, 38(3):1287--1319.

\bibitem[Roos, 1999]{roos1999rate}
Roos, B. (1999).
\newblock On the rate of multivariate {P}oisson convergence.
\newblock {\em Journal of Multivariate Analysis}, 69(1):120--134.

\bibitem[Rothschild and Stiglitz, 1978]{rothschild1978increasing}
Rothschild, M. and Stiglitz, J.~E. (1978).
\newblock Increasing risk: I. a definition.
\newblock In {\em Uncertainty in Economics}, pages 99--121. Elsevier.

\bibitem[Roverato et~al., 2013]{roverato2013log}
Roverato, A., Lupparelli, M., and La~Rocca, L. (2013).
\newblock Log-mean linear models for binary data.
\newblock {\em Biometrika}, 100(2):485--494.

\bibitem[Saoub, 2021]{saoub2021graph}
Saoub, K.~R. (2021).
\newblock {\em Graph Theory: An Introduction to Proofs, Algorithms, and
  Applications}.
\newblock CRC Press.

\bibitem[Shaked and Shanthikumar, 2007]{shaked2007}
Shaked, M. and Shanthikumar, J.~G. (2007).
\newblock {\em Stochastic Orders}.
\newblock Springer.

\bibitem[Steutel and van Harn, 1979]{steutel1979discrete}
Steutel, F.~W. and van Harn, K. (1979).
\newblock Discrete analogues of self-decomposability and stability.
\newblock {\em Annals of Probability}, 7(5):893--899.

\bibitem[Tasche, 2007]{tasche2007capital}
Tasche, D. (2007).
\newblock Capital allocation to business units and sub-portfolios: the {E}uler
  principle.
\newblock {\em arXiv preprint arXiv:0708.2542}.

\bibitem[Teugels, 1990]{teugels1990some}
Teugels, J.~L. (1990).
\newblock Some representations of the multivariate {B}ernoulli and binomial
  distributions.
\newblock {\em Journal of Multivariate Analysis}, 32(2):256--268.

\bibitem[Theodoridis, 2015]{theodoridis2015machine}
Theodoridis, S. (2015).
\newblock {\em Machine Learning: A Bayesian and Optimization Perspective}.
\newblock Academic Press.

\bibitem[Venema, 1992]{venema1992modeling}
Venema, H.~W. (1992).
\newblock Modeling fiber type grouping by a binary {M}arkov random field.
\newblock {\em Muscle \& Nerve: Official Journal of the American Association of
  Electrodiagnostic Medicine}, 15(6):725--732.

\bibitem[Wainwright and Jordan, 2008]{wainwright2008graphical}
Wainwright, M.~J. and Jordan, M.~I. (2008).
\newblock Graphical models, exponential families, and variational inference.
\newblock {\em Foundations and Trends{\textregistered} in Machine Learning},
  1(1--2):1--305.

\bibitem[Wakeland, 2017]{wakeland2017exploring}
Wakeland, K.~W. (2017).
\newblock {\em Exploring Dependence in Binary Markov Random Field Models}.
\newblock PhD thesis, Iowa State University.

\end{thebibliography}

	\appendix

	\section{Proofs}

    \subsection{Proof of Proposition~\ref{prop:corr}}
    \label{sect:proofCorr}
    The case where $(u,v)\in\mathcal{E}$ is trivial due to the mean parameterization. We prove for the general case where $k\geq 3$ vertices are on $\mathrm{path}(u,v)$. Let us denote the vertices on the path by $\{\ell_i, \; i \in \{1,\ldots,k\}\}$ such that $\ell_1=u$ and $\ell_k=v$. Hence, $\mathrm{Cov}(J_u,J_v) = \mathrm{Cov}(J_{\ell_1}, J_{\ell_k})$. By conditioning on $J_{\ell_{k-1}}$, 
	\begin{equation}
		\mathrm{Cov}(J_{\ell_1}, J_{\ell_k}) = \mathrm{Cov}(\mathrm{E}[J_{\ell_1}|J_{\ell_{k\text{-}1}}],\mathrm{E}[J_{\ell_k}|J_{\ell_{k\text{-}1}}]) + \mathrm{E}[\mathrm{Cov}(J_{\ell_1}, J_{\ell_k}|J_{\ell_{k\text{-}1}})].
		\label{eq:CovNj1Nj2Noedge-2}
	\end{equation}
	The second term in (\ref{eq:CovNj1Nj2Noedge-2}) is equal to zero by the global Markov property. Note that, from the mean parameterization fixing bivariate distributions associated to edges, for any $m\in\{2,\ldots, k\}$,
    \begin{align*}
        \mathrm{E}[J_{\ell_m}|J_{\ell_{m\text{-}1}}] &= \left(1-J_{\ell_{m\text{-}1}}\right)\left(q_{\ell_m} - \frac{\alpha_{(\ell_{m\text{-}1},\ell_m)}\sigma_{(\ell_{m\text{-}1},\ell_m)}}{1-q_{\ell_{m\text{-}1}}} \right) + J_{\ell_{m\text{-}1}} \left(q_{\ell_m}+\frac{\alpha_{(\ell_{m\text{-}1},\ell_m)}\sigma_{(\ell_{m\text{-}1},\ell_m)}}{q_{\ell_{m\text{-}1}}}  \right)\\ &= q_{\ell_m} - \frac{\alpha_{(\ell_{m\text{-}1},\ell_m)}\sigma_{(\ell_{m\text{-}1},\ell_m)}}{1-q_{\ell_{m\text{-}1}}} + \frac{\alpha_{(\ell_{m\text{-}1},\ell_m)}\sigma_{(\ell_{m\text{-}1},\ell_m)}}{q_{\ell_{m\text{-}1}}(1-q_{\ell_{m\text{-}1})}}J_{\ell_{m\text{-}1}}. 
    \end{align*}
    Therefore, \eqref{eq:CovNj1Nj2Noedge-2} becomes
    \begin{align*}
        \mathrm{Cov}(J_{\ell_1}, J_{\ell_k}) &= \mathrm{Cov}\left(\mathrm{E}[J_{\ell_1}|J_{\ell_{k\text{-}1}}], q_{\ell_m} - \frac{\alpha_{(\ell_{k\text{-}1},\ell_k)}\sigma_{(\ell_{k\text{-}1},\ell_k)}}{1-q_{\ell_{k\text{-}1}}} + \frac{\alpha_{(\ell_{k\text{-}1},\ell_k)}\sigma_{(\ell_{k\text{-}1},\ell_k)}}{q_{\ell_{k\text{-}1}}(1-q_{\ell_{k\text{-}1})}}J_{\ell_{k\text{-}1}} \right)\\
        &= \frac{\alpha_{(\ell_{k\text{-}1},\ell_k)}\sigma_{(\ell_{k\text{-}1},\ell_k)}}{q_{\ell_{k\text{-}1}}(1-q_{\ell_{k\text{-}1})}}\mathrm{Cov}\left(\mathrm{E}[J_{\ell_1}|J_{\ell_{k\text{-}1}}], J_{\ell_{k\text{-}1}} \right).
    \end{align*}
    Repeating this rationale iteratively, we obtain 
    \begin{align}
        \mathrm{Cov}(J_{\ell_1}, J_{\ell_k}) = \left(\prod_{m=2}^k\frac{\alpha_{(\ell_{m\text{-}1},\ell_m)}\sigma_{(\ell_{m\text{-}1},\ell_m)}}{q_{\ell_{m\text{-}1}}(1-q_{\ell_{m\text{-}1})}}\right)\mathrm{Var}(J_{\ell_1}) &=  \left(\prod_{m=2}^k \alpha_{(\ell_{m\text{-}1},\ell_m)} \right) \left(\frac{\sigma_{(\prod_{m=2}^k \ell_{m\text{-}1},\ell_m)}}{\prod_{j=2}^{m-1} q_{\ell_{j}(1-q_{\ell_{j\text{+}1})}}}\right)\notag\\ &= \left(\prod_{m=2}^k \alpha_{(\ell_{m\text{-}1},\ell_m)} \right) \sqrt{q_{l_1}(1-q_{l_1})q_{l_k}(1-q_{l_k})}.
        \label{eq:lllllll}
    \end{align}
    Reverting back to the original indices, \eqref{eq:lllllll} is rewritten as
    \begin{equation*}
        \mathrm{Cov}(J_{u}, J_{v}) = \left(\prod_{e\in \mathrm{path}(u,v)} \alpha_{e} \right) \sqrt{q_{u}(1-q_{u})q_{v}(1-q_{v})}. 
    \end{equation*}
     The desired result directly ensues. \hfill$\square$

		\subsection{Proof of Theorem \ref{th:PoissonApprox}}	
	\label{sect:proofPoissonApprox}
The proof follows a Lindeberg-type replacement scheme in which the Poisson model is transformed into the Bernoulli model one vertex at a time.   
Choose any root $r\in\mathcal{V}$ and define all filiation relations with respect to this root. For every $v\in\mathrm{dsc}(r)$, set $N_v^{\prime}=
\alpha_{(\mathrm{pa}(v),v)}\circ J_{\mathrm{pa}(v)}+L_v$. This auxiliary variable combines the Bernoulli parent from the Ising model with the Poisson innovation from the MPMRF. It is introduced to replace the Poisson components by Bernoulli components successively, while preserving a coherent and well defined MRF throughout the proof. For ease of exposition, assume that the vertices are labeled in increasing order away from the root, so that $\mathrm{pa}(v) < v$ for every $v\in\mathcal{V}$.

\bigskip

Let $\boldsymbol{N}^{(0)} = \boldsymbol{N}$  and, for each $i\in\{1,\ldots, d\}$, let $\boldsymbol{N}^{(i)}$ be obtained from $\boldsymbol{N}^{(i-1)}$ by changing its $i$th component to $J_i$; this transforms $N_j$ to $N^{\prime}_j$ for all $j\in\mathrm{ch}(i)$. Also, note that $\boldsymbol{N}^{(d)} = \boldsymbol{J}$. 
By a triangle inequality, 
\begin{equation}
    \mathbbm{d}_{\mathrm{TV}}(\boldsymbol{N}, \boldsymbol{J}) \leq \sum_{i=1}^d \mathbbm{d}_{\mathrm{TV}}\left(\boldsymbol{N}^{(i-1)}, \boldsymbol{N}^{(i)}\right).
\label{eq:triangle-inequality}
\end{equation}

    For a set $\Omega$, let $\mathscr{P}(\Omega)$ denote its power set. 
		Because the support of $\boldsymbol{N}^{(i)}$ is $\{0,1\}^i\times\mathbb{N}^{d-i}\subseteq \mathrm{supp}\boldsymbol{N}^{(i-1)} = \{0,1\}^{i-1}\times \mathbb{N}^{d-i+1}$, we have
		\begin{align}
			\mathbbm{d}_{\mathrm{TV}}\left(\boldsymbol{N}^{(i-1)}, \boldsymbol{N}^{(i)}\right) &= \sup_{A\in\mathscr{P}
            (\{0,1\}^{i-1}\times \mathbb{N}^{d-i+1})}\left\{ \mathrm{Pr}\left(\boldsymbol{N}^{(i-1)}\in A\right) - \mathrm{Pr}\left(\boldsymbol{N}^{(i)}\in A\right) \right\} \notag\\
            &= \sup_{A\in\mathscr{P}
            (\{0,1\}^i\times\mathbb{N}^{d-i})}\left\{ \mathrm{Pr}\left(\boldsymbol{N}^{(i-1)}\in A\right) - \mathrm{Pr}\left(\boldsymbol{N}^{(i)}\in A\right) \right\}  + \mathrm{Pr}\left(N_i^{(i-1)} \geq 2\right) 
			\label{eq:proofPoissonApprox-Goal}
		\end{align}
for every $i\in\{1,\ldots,d\}$. For the first term in \eqref{eq:proofPoissonApprox-Goal}, the supremum is taken over events contained in $\{0,1\}^{i}\times\mathbb N^{d-i}$, that is, over events for which the $i$th coordinate can only take the values $0$ and $1$. On this restricted state space, the vectors $\boldsymbol N^{(i-1)}$ and $\boldsymbol N^{(i)}$ have identical stochastic constructions except for their $i$th components, namely $N_i^{(i-1)}$ and $N_i^{(i)}$. We partitioned the state space according to whether the $i$th component exceeds 1. Note that $\mathrm{Pr}(N_i^{(i)}\geq 2)=0$, and thus the part exceeding 1 yields $\mathrm{Pr}\left(N_i^{(i-1)} \geq 2\right)$. The part less or equal to 1 is what is accounted by the first supremum of \eqref{eq:proofPoissonApprox-Goal}. Therefore, 
\begin{align*}
    \sup_{A\in\mathscr{P}
            (\{0,1\}^i\times\mathbb{N}^{d-i})}\left\{ \mathrm{Pr}\left(\boldsymbol{N}^{(i-1)}\in A\right) - \mathrm{Pr}\left(\boldsymbol{N}^{(i)}\in A\right) \right\} = \sup_{A\in\mathscr{P}(\{0,1\})}\left\{ \mathrm{Pr}\left(N_i^{(i-1)}\in A\right) - \mathrm{Pr}\left(N_i^{(i)}\in A\right) \right\}
\end{align*}
For $i=r=1$, 
\begin{align}
    &\sup_{A\in\mathscr{P}(\{0,1\})}\left\{ \mathrm{Pr}\left(N_1^{(0)}\in A\right) - \mathrm{Pr}\left(N_1^{(1)}\in A\right) \right\} =  \sup_{A\in\mathscr{P}(\{0,1\})}\left\{ \mathrm{Pr}\left(N_1\in A\right) - \mathrm{Pr}\left(J_1\in A\right) \right\}\notag\\
    &\quad\quad= \mathrm{Pr}\left(N_1 = 0\right) - \mathrm{Pr}\left(J_1 = 0\right)
    = \mathrm{e}^{-q} - (1-q). \label{eq:PoissonApprox-1-r}
\end{align}
For $i>1$, 
\begin{align}
    &\sup_{A\in\mathscr{P}(\{0,1\})}\left\{ \mathrm{Pr}\left(N_i^{(i-1)}\in A\right) - \mathrm{Pr}\left(N_i^{(i)}\in A\right) \right\} =  \sup_{A\in\mathscr{P}(\{0,1\})}\left\{ \mathrm{Pr}\left(N_i^{\prime}\in A\right) - \mathrm{Pr}\left(J_i\in A\right) \right\}\notag\\
    &\quad\quad= \mathrm{Pr}\left(N_i^{\prime} = 0\right) - \mathrm{Pr}\left(J_i = 0\right)\notag\\
    &\quad \quad= (1-q + q(1-\alpha_{(\mathrm{pa}(i),i)}))\mathrm{e}^{-(1-\alpha_{(\mathrm{pa}(i),i)})q} - (1-q). \label{eq:PoissonApprox-1-i}
\end{align}
In both \eqref{eq:PoissonApprox-1-r} and \eqref{eq:PoissonApprox-1-i}, the second equality follows by exhaustion of all four possible combinations. 
 We show that the right-hand side of \eqref{eq:PoissonApprox-1-i} reaches is dominated by the value obtained for the case $i=1$, that is
\begin{equation}
    (1-q + q(1-\alpha_{(\mathrm{pa}(i),i)}))\mathrm{e}^{-(1-\alpha_{(\mathrm{pa}(i),i)})q} - (1-q) \leq \mathrm{e}^{-q} - (1-q).
\label{eq:obj-beta}
\end{equation}
Let $\beta=1-\alpha_{(\mathrm{pa}(i),i)}$, thus $\beta\in[0,1]$, and define the right-hand side of \eqref{eq:PoissonApprox-1-i} as a function $\varphi$ of $\beta$, namely
\begin{equation*}
    \varphi(\beta)= (1-q + q \beta )\mathrm{e}^{-\beta q} - (1-q).
\end{equation*}
Given that
\begin{equation*}
    \frac{\mathrm{d}}{\mathrm{d}\beta}\varphi(\beta)=q\mathrm{e}^{-\beta q}-q(1-q+q \beta)\mathrm{e}^{-\beta q} 
    = q^2 \mathrm{e}^{-\beta q}(1-\beta)\geq 0,
\end{equation*}
the function $\phi$ is increasing in $\beta$, with $\beta\in[0,1]$. It thus reaches its maximum value at $\beta=1$, yielding \eqref{eq:obj-beta}. 

\bigskip

For the second term in \eqref{eq:proofPoissonApprox-Goal}, we have, for $i=r=1$, 
\begin{equation}
    \mathrm{Pr}\left(N_i^{(i-1)}\geq 2\right) = \mathrm{Pr}\left(N_1\geq 2\right) = 1-\mathrm{e}^{-q} - q\mathrm{e}^{-q} , \label{eq:PoissonApprox-2-1} 
\end{equation}
and, for $i>1$, 
\begin{align}
&\mathrm{Pr}\left(N_i^{(i-1)}\geq 2\right) = 
    1 - \mathrm{Pr}\left(N_i^{\prime}=0\right) - \mathrm{Pr}\left(N_i^{\prime}=1\right)\notag\\
    & = 1 - (1-q + q(1-\alpha_{(\mathrm{pa}(i),i)}) + (1-q)q(1-\alpha_{(\mathrm{pa}(i),i)}) + q^2(1-\alpha_{(\mathrm{pa}(i),i)})^2 + q\alpha_{(\mathrm{pa}(i),i)})\mathrm{e}^{-(1-\alpha_{(\mathrm{pa}(i),i)})q} \notag\\
    \label{eq:PoissonApprox-2-i}
\end{align}
We compare the probability $\mathrm{Pr}(N_i^{(i-1)}\geq 2)$ in a similar way as for the first term in \eqref{eq:proofPoissonApprox-Goal}. First, let $\beta=1-\alpha_{(\mathrm{pa}(i),i)}$ in \eqref{eq:PoissonApprox-2-i}. Let the right-hand side of \eqref{eq:PoissonApprox-2-i} be expressed as a function $\psi$ of $\beta$, for $\beta\in[0,1]$. We have
\begin{align*}
\psi(\beta)&= 
     1 - \{1-q + q \beta + (1-q)q \beta + q^2 \beta^2 + q(1-\beta) \}\mathrm{e}^{-\beta q} \\
     &=1 - \{1+(1-q)q\beta + q^2 \beta^2  \}\mathrm{e}^{-\beta q}.
\end{align*}
Note that
\begin{align*}
    \psi^{\prime}(\beta)&=(-(1-q)q-2q^2\beta)\mathrm{e}^{-\beta q}+\{(1+(1-q)q\beta+q^2 \beta^2)\}q \mathrm{e}^{-\beta q} \\
    &=\{q^2 (1-\beta)-q^3\beta (1-\beta)\}\mathrm{e}^{-\beta q} =(q^2(1-\beta)(1-q))\mathrm{e}^{-\beta q} \geq 0,
\end{align*}
and hence, the function $\psi$ is increasing in $\beta$. Its maximum value is at $\beta=1$ and it is equal to
\begin{equation*}
    \psi(1)=1 - \{1+(1-q)q + q^2   \}\mathrm{e}^{-q},
\end{equation*}
which corresponds to the value of $\mathrm{Pr}(N_i^{(i-1)}\geq 2)$ for $i=1$.

Combining this with \eqref{eq:triangle-inequality}, \eqref{eq:proofPoissonApprox-Goal}, \eqref{eq:PoissonApprox-1-r}, \eqref{eq:PoissonApprox-1-i}, \eqref{eq:obj-beta} and \eqref{eq:PoissonApprox-2-1} leads to
		\begin{equation*}
			\mathbbm{d}_{\mathrm{TV}}(\boldsymbol{N}, \boldsymbol{J}) \leq \sum_{i=1}^d (\mathrm{e}^{-q}-(1-q) + 1-(1+q)\mathrm{e}^{-q}) = dq(1-\mathrm{e}^{-q}).
		\end{equation*}

		That is, as $q$ approaches 0, the total variation distance between $\boldsymbol{J}$ and $\boldsymbol{N}$ abates at the same rate as $q(1-\mathrm{e}^{-q})$. 
        \hfill $\square$

	\subsection{Proof of Theorem \ref{th:PoissonApproxConvexOrder}}	
	\label{sect:proofPoissonApproxConvexOrder}

     Similarly to the previous proof, choose any root $r\in\mathcal{V}$ and, for every $v\in\mathrm{dsc}(r)$ define an auxiliary random variable $N_v^{\prime}=\alpha_{(\mathrm{pa}(v),v)}\circ J_{\mathrm{pa}(v)} + L_v$ which combines the Bernoulli parent with the Poisson offspring mechanism. This hybrid variable enables us to replace the vertices one at a time while preserving the dependence structure required; that is the MRF remains well-defined.  
        First, we prove $(J_{v} | J_{\mathrm{pa}(v)} = x_{\mathrm{pa}(v)}) \preceq_{cx} (N_{v}^{\prime} | J_{\mathrm{pa}(v)} = x_{\mathrm{pa}(v)})$, $v \in \mathrm{dsc}(r)$, $x_{\mathrm{pa}(v)}\in\{0,1\}$. We have 
		\begin{align}
			&\mathrm{Pr}(J_v=0|J_{\mathrm{pa}(v)} = 0) = 1-(1-\alpha_{(\mathrm{pa}(v),v)})q \leq \mathrm{e}^{-(1-\alpha_{(\mathrm{pa}(v),v)})q} = \mathrm{Pr}(N_v^{\prime}=0|J_{\mathrm{pa}(v)}=0)\notag\\
			&\mathrm{Pr}(J_v\leq 1|J_{\mathrm{pa}(v)} = 0) = 1 \geq (1+(1-\alpha_{(\mathrm{pa}(v),v)})q)\mathrm{e}^{-(1-\alpha_{(\mathrm{pa}(v),v)})q} = \mathrm{Pr}(N_v^{\prime}\leq 1|J_{\mathrm{pa}(v)}=0),\notag
		\end{align}
		which implies $(J_v|J_{\mathrm{pa}(v)}=0)\preceq_{cx}(N_v^{\prime}|J_{\mathrm{pa}(v)}=0)$ by the Karlin-Novikoff sufficient condition for convex ordering (see e.g.~Theorem~1.5.17 of \cite{muller2002}) and the fact that both sides of the order relation have the same mean. Similarly, 
		\begin{align}
			&\mathrm{Pr}(J_v=0|J_{\mathrm{pa}(v)} = 1) = (1-\alpha_{(\mathrm{pa}(v),v)})(1-q) \leq (1-\alpha_{(\mathrm{pa}(v),v)}) \mathrm{e}^{-(1-\alpha_{(\mathrm{pa}(v),v)})q} = \mathrm{Pr}(N_v^{\prime}=0|J_{\mathrm{pa}(v)}=1)\notag\\
			&\mathrm{Pr}(J_v\leq 1|J_{\mathrm{pa}(v)} = 1) = 1 \geq (1 +(1-\alpha_{(\mathrm{pa}(v),v)})^2q)\mathrm{e}^{-(1-\alpha_{(\mathrm{pa}(v),v)})q} = \mathrm{Pr}(N_v^{\prime}\leq 1|J_{\mathrm{pa}(v)}=1),\notag
		\end{align}
		which implies $(J_v|J_{\mathrm{pa}(v)}=1)\preceq_{cx}(N_v^{\prime}|J_{\mathrm{pa}(v)}=1)$ again by Karlin-Novikoff and their having the same mean. Combining both results renders as
		\begin{equation}
			(J_v|J_{\mathrm{pa}(v)}=x_{\mathrm{pa(v)}})\preceq_{cx}(N_v^{\prime}|J_{\mathrm{pa}(v)}=x_{\mathrm{pa}(v)}),\quad \quad x_{\mathrm{pa}(v)}\in\{0,1\}.
			\label{eq:KcxM-JcxNprime}
		\end{equation}
		Consequently, 		
		\begin{equation}
			((J_{\mathrm{pa}(v)} + J_v)|J_{\mathrm{pa}(v)}=x_{\mathrm{pa}(v)})\preceq_{cx}((J_{\mathrm{pa}(v)} + N_v^{\prime})|J_{\mathrm{pa}(v)}=x_{\mathrm{pa}(v)}),\quad \quad x_{\mathrm{pa}(v)}\in\{0,1\}.
			\label{eq:KcxM-ExJpa+J}
		\end{equation}
        The next part of the proof is predicated on the notion of stochastic increasingness and convexity. We show that random variables following the construction of either models (or a hybrid of the two) satisfy this property, and rely on this to add random variables to both sides of \eqref{eq:KcxM-ExJpa+J} while preserving the convex order.  
		A random variable $X$ is stochastically increasing and convex (SICX) with respect to a random variable $\Theta$ if $\mathrm{E}[\varphi(X)|\Theta=\theta]$ is increasing convex in $\theta$ for every increasing convex function $\varphi$. 
       Since binomial thinning satisfies the semigroup property with respect to the number of trials, the thinning operator preserves SICX (about the semigroup property and SICX, see Example 8.A.7 of \cite{shaked2007}). 
        For notation purposes, let $\Phi_v = \alpha_{(\mathrm{pa}(v),v)}\circ\Theta_{\mathrm{pa}(v)} + L_v$, with $\Theta_{\mathrm{pa}(v)}$ acting as either $J_{\mathrm{pa}(v)}$, $N_{\mathrm{pa}(v)}^{\prime}$ or $N_{\mathrm{pa}(v)}$; we use $\Phi_v$ to consider the three cases simultaneously. From the closure of SICX-ness on sums of conditionally independent random variables, as stated by Theorem 8.A.15 of \cite{shaked2007},  $(\Theta_{\mathrm{pa}(v)} + \alpha_{(\mathrm{pa}(v),v)}\circ\Theta_{\mathrm{pa}(v)}+L_v)$, \textit{i.e.} $(\Theta_{\mathrm{pa}(v)} + \Phi_v)$, is SICX with respect to $ \Theta_{\mathrm{pa}(v)}$, which, from Theorem 8.A.17 of \cite{shaked2007}, the closure of SICX-ness on mixtures, implies $\left(\Theta_{\mathrm{pa}(v)}+\Phi_v + \sum_{j\in\mathrm{dsc}(v)} N_j\right)$ is SICX with respect to $\Theta_{\mathrm{pa}(v)}$. Invoking again the closure of SICX-ness on sums of conditionally independent random variables, we deduce 
		\begin{equation}
			\left(\Theta_{u}+\sum_{i\in\mathrm{ch}(u)}\left(\Phi_i + \sum_{j\in\mathrm{dsc}(i)} N_j\right)\right) \quad \text{is SICX with respect to}\;\Theta_{u},\, u\in\mathcal{V}. 
			\label{eq:KcxM-SICXNstar}
		\end{equation} 
		 Applying Theorem 8.A.14a of \cite{shaked2007} with $(J_v\mid J_{\mathrm{pa}(v)} = x_{\mathrm{pa}(v)})$ and $(N^{\prime}_v \mid J_{\mathrm{pa}(v)} = x_{\mathrm{pa}(v)})$ as the mixing random variables, we have, given (\ref{eq:KcxM-JcxNprime}), (\ref{eq:KcxM-SICXNstar}), and that both sides share the same expectation, for any $v\in\mathrm{dsc}(r)$,
		{\setlength{\abovedisplayskip}{3pt}
			\begin{equation}
				{\scalebox{0.80}{
						\begin{minipage}{\linewidth}
							\begin{equation*}
								\left(\left.\left(J_v + \sum_{i\in\mathrm{ch}(v)}\left(N_i^{\prime} + \sum_{j\in\mathrm{dsc}(i)} N_j\right)\right)\right|J_{\mathrm{pa}(v)}=x_\mathrm{pa(v)}\right)\preceq_{cx} \left(\left.\left(N_v^{\prime} + \sum_{i\in\mathrm{ch}(v)}\left(N_i + \sum_{j\in\mathrm{dsc}(i)} N_j\right)\right)\right|J_{\mathrm{pa}(v)}=x_\mathrm{pa(v)}\right), 
				\end{equation*}\end{minipage}}}\quad\quad\quad\label{eq:KcxM-RelationImportante}
		\end{equation}}
		
		for $x_{\mathrm{pa}(v)}\in\{0,1\}$. Let $\psi=\max_{v\in\mathcal{V}}|\mathrm{path}(r,v)|$. Partition the vertices of the underlying tree in subsets $\{\Psi_k,\; k\in\{0,\ldots,\psi\}\}$ of $\mathcal{V}$, depending on the length of their path to the root. That is $j\in\Psi_i$ if $|\mathrm{path}(r,j)|=i$. Then, in the same spirit as in (\ref{eq:KcxM-ExJpa+J}) and given both models satisfy the global Markov property, we have, given (\ref{eq:KcxM-JcxNprime}),
		\begin{equation}  
			\left(\left.\left(\sum_{i=0}^{\psi}\sum_{j\in\Psi_i}J_j  \right)\right|\bigcap_{i=0}^{\psi-1}\bigcap_{j\in\Psi_i}\{J_j = x_j\}\right) \preceq_{cx} \left( \left.\left(\sum_{i=0}^{\psi-1}\sum_{j\in\Psi_i}J_j + \sum_{k\in\Psi_{\psi}}N_k^{\prime} \right)\right|\bigcap_{i=0}^{\psi-1}\bigcap_{j\in\Psi_i}\{J_j = x_j\}\right),
			\label{eq:KcxM-Psi1}
		\end{equation}
		for $\boldsymbol{x}\in\{0,1\}^d$. By multiple applications of Theorem 3.A.12b of \cite{shaked2007}, successively using every $J_j$, $j\in\Psi_{\psi-1}$, to act as the mixture random variable, (\ref{eq:KcxM-Psi1}) becomes \\
		{\scalebox{0.95}{
				\begin{minipage}{\linewidth}
					\begin{align*}  
						\left(\left.\left(\sum_{i=0}^{\psi}\sum_{j\in\Psi_i}J_j  \right)\right|\bigcap_{i=0}^{\psi-2}\bigcap_{j\in\Psi_i}\{J_j = x_j\}\right) &\preceq_{cx}\left( \left.\left(\sum_{i=0}^{\psi-1}\sum_{j\in\Psi_i}J_j + \sum_{k\in\Psi_{\psi}}N_k^{\prime} \right)\right|\bigcap_{i=0}^{\psi-2}\bigcap_{j\in\Psi_i}\{J_j = x_j\}\right)\notag\\
						&\preceq_{cx} \left(\left.\left(\sum_{i=0}^{\psi-2}\sum_{j\in\Psi_i}J_j +  \sum_{k\in\Psi_{\psi-1}}N_k^{\prime} + \sum_{l\in\Psi_{\psi}}N_l \right)\right|\bigcap_{i=0}^{\psi-2}\bigcap_{j\in\Psi_i}\{J_j = x_j\}\right),
		\end{align*}\end{minipage}}}
		given (\ref{eq:KcxM-RelationImportante}). One then repeats the multiple applications of Theorem 3.A.12b, this time using every $J_j$, $j\in\Psi_{\psi-2}$, and obtains, for $\boldsymbol{x}\in\{0,1\}^d$,\\
		{\scalebox{0.85}{
				\begin{minipage}{\linewidth}
					\begin{equation*}
						\left(\left.\left(\sum_{i=0}^{\psi}\sum_{j\in\Psi_i}J_j  \right)\right|\bigcap_{i=0}^{\psi-3}\bigcap_{j\in\Psi_i}\{J_j = x_j\} \right)
						\preceq_{cx} \left(\left.\left(\sum_{i=0}^{\psi-2}\sum_{j\in\Psi_i}J_j +  \sum_{k\in\Psi_{\psi-1}}N_k^{\prime} + \sum_{l\in\Psi_{\psi}}N_l \right)\right|\bigcap_{i=0}^{\psi-3}\bigcap_{j\in\Psi_i}\{J_j = x_j\}\right) 
					\end{equation*}
		\end{minipage}}}
		{\setlength{\abovedisplayskip}{3pt}
			\begin{equation}
				{\scalebox{0.85}{
						\begin{minipage}{\linewidth}
							\begin{equation*}
								\mathcolor{White}{	\left(\left.\left(\sum_{i=0}^{\psi}\sum_{j\in\Psi_i}J_j  \right)\right|\bigcap_{i=0}^{\psi-3}\bigcap_{j\in\Psi_i}\{ = x_j\} \right)}
								\preceq_{cx}
								\left(\left.\left(\sum_{i=0}^{\psi-3}\sum_{j\in\Psi_i}J_j +  \sum_{k\in\Psi_{\psi-2}}N_k^{\prime} + \sum_{r=\psi-1}^\psi\sum_{l\in\Psi_{r}}N_l \right)\right|\bigcap_{i=0}^{\psi-3}\bigcap_{j\in\Psi_i}\{J_j = x_j\}\right),
				\end{equation*}\end{minipage}}}\quad\quad\quad	\label{eq:KcxM-Psi3}
		\end{equation}}
		again given (\ref{eq:KcxM-RelationImportante}). We repeat this argument $\psi-3$ other times, such that (\ref{eq:KcxM-Psi3}) unravels to
		\begin{align*}
			\left(\left.\left(\sum_{i=0}^{\psi}\sum_{j\in\Psi_i}J_j  \right)\right|\bigcap_{j\in\Psi_0}\{J_j = x_j\} \right)
			&\preceq_{cx} \left(\left.\left(\sum_{j\in\Psi_0}J_j +  \sum_{k\in\Psi_{1}}N_k^{\prime} + \sum_{r=2}^\psi\sum_{l\in\Psi_{r}}N_l \right)\right|\bigcap_{j\in\Psi_0}\{J_j = x_j\}\right).
		\end{align*}
		Note that the sole vertex in subset $\Psi_0$ is the root. Applying one last time Theorem 3.A.12b and dropping the $\Psi_k$ notation, we obtain
		\begin{align}
			\left(\sum_{v\in\mathcal{V}} J_v\right) \preceq_{cx} \left(J_r + \sum_{j\in\mathrm{ch(r)}}N_j^{\prime} + \sum_{v\in\mathcal{V}\backslash(\{r\}\cup\mathrm{ch}(r))} N_v\right). 
			\label{eq:KcxM-AlmostDone}
		\end{align}
		As exposed by (\ref{eq:KcxM-SICXNstar}), $\left(\Theta_r + \sum_{j\in\mathrm{ch(r)}}\Phi_j + \sum_{v\in\mathcal{V}\backslash(\{r\}\cup\mathrm{ch}(r))} N_v\right)$ is SICX with respect to $\Theta_r$. Since $J_r\preceq_{cx} N_r$, since their respective distributions are Bernoulli and Poisson with identical means, applying Theorem 8.A.14a of \cite{shaked2007} yields
		\begin{align}
			\left(J_r + \sum_{j\in\mathrm{ch(r)}}N_j^{\prime} + \sum_{v\in\mathcal{V}\backslash(\{r\}\cup\mathrm{ch}(r))} N_v\right) \preceq_{cx} \left(\sum_{v\in\mathcal{V}} N_v\right).
			\label{eq:KcxM-Done}
		\end{align}
		Joining (\ref{eq:KcxM-AlmostDone}) and (\ref{eq:KcxM-Done}), we deduce $\left(\sum_{v\in\mathcal{V}}J_v\right)\preceq_{cx}\left(\sum_{v\in\mathcal{V}}N_v\right)$. This is $K\preceq_{cx}M$.
	\hfill $\square$

\end{document}